\documentclass[11pt]{article}
\usepackage{amsfonts}
\usepackage{amsfonts,latexsym,amsmath,amscd,geometry}
\usepackage{amssymb}
\usepackage{latexsym}

\newcommand \nc{\newcommand}
\newtheorem{theorem}{Theorem}[section]
\newtheorem{lemma}[theorem]{Lemma}
\newtheorem{proposition}[theorem]{Proposition}
\newtheorem{corollary}[theorem]{Corollary}

\newtheorem{remark}[theorem]{Remark}

\nc{\ba}{\begin{array}}\nc{\ea}{\end{array}}
\nc{\be}{\begin{eqnarray}}\nc{\ee}{\end{eqnarray}}
\nc{\beq}{\begin{equation}}\nc{\eeq}{\end{equation}}
\nc{\bex}{\begin{eqnarray*}}\nc{\eex}{\end{eqnarray*}}
\nc{\btm}{\begin{theorem}} \nc{\etm}{\end{theorem}}
\nc{\blm}{\begin{lemma}} \nc{\elm}{\end{lemma}}
\nc{\R}{\mathbb{R}} \nc{\va}{\varepsilon} \nc{\ls}{\limits}

\def\pf{\noindent{\bf Proof.\quad}}\def\endpf{\hfill$\Box$}

\begin{document}
\title{{\bf Global well-posedness and zero-diffusion limit of classical solutions to the 3D conservation
 laws arising in chemotaxis}}
\author{Hongyun Peng\thanks{School
of Mathematics and Statistics, Central China Normal University,
Wuhan 430079, China. Email: penghy010@163.com},\quad Huanyao
Wen\thanks{School of Mathematics and Statistics, Central China
Normal University, Wuhan 430079, China, and Department of Petroleum
Engineering, Faculty of Science and Technology, University of
Stavanger, 4036 Stavanger, Norway. Email:
huanyaowen@hotmail.com},\quad Changjiang Zhu\thanks{Corresponding
author. School of Mathematics and Statistics, Central China Normal
University, Wuhan 430079, China. Email: cjzhu@mail.ccnu.edu.cn}  }
\date{}

\maketitle

\begin{abstract}
In this paper, we study the relationship between a diffusive model
and a non-diffusive model which are both derived from the well-known
Keller-Segel model, as a coefficient of diffusion $\varepsilon$ goes
to zero. First, we establish the global well-posedness of classical
solutions to the Cauchy problem for the diffusive model
 with smooth initial
data which is of small $L^2$ norm, together with some {\it a priori}
estimates uniform for $t$ and $\varepsilon$. Then we investigate the
zero-diffusion limit, and get the global well-posedness of classical
solutions to the Cauchy problem for the non-diffusive model.
Finally, we derive the convergence rate of the diffusive model
toward the non-diffusive model. It is shown that the convergence
rate in $L^\infty$ norm is of the order
$O\left(\varepsilon^{\frac{1}{2}}\right)$. It should be noted that
the initial data is small in $L^2$-norm but can be of large
oscillations with constant state at far field. As a byproduct, we
improve the corresponding result on the well-posedness of the
non-difussive model which requires small oscillations.
\end{abstract}


\vspace{4mm}
\section {Introduction}
\setcounter{equation}{0} \setcounter{theorem}{0} In this paper, we
are interested in a system of conservation laws arising in
chemotaxis
\begin{equation}\label{Equation@}
\left\{
\begin{array}{l}
p^\varepsilon_t-\nabla\cdot\left(p^\varepsilon
\bf{q}^\varepsilon\right)
=D\triangle p^\varepsilon,\\[3mm]
{\bf{q}}^\varepsilon_t +\nabla\left(\varepsilon
|{\bf{q}}^\varepsilon|^2-p^\varepsilon\right)=\varepsilon \triangle
{\bf{q}}^\varepsilon,
\end{array}
\right.
\end{equation}
with initial data
\begin{equation}\label{Initial}
\left(p^\varepsilon, {\bf
q}^\varepsilon\right)(x,0)=\left(p_0(x),{\bf
q}_0(x)\right)\rightarrow(p_\infty, 0)\ \ {\rm as}\ \ |x|\rightarrow
\infty.
\end{equation}

The chemotaxis model was preoposed by Keller and Segel in
\cite{Keller71} to describe the traveling band behavior of bacteria
due to the chemotactic response observed in experiments
\cite{Adler66,Adler69}. The following Keller-Segel model has been
extensively studied
\begin{equation}\label{Keller-Segel1}
\left\{
\begin{array}{l}
u_t=\nabla\cdot \left(D\nabla u-\chi u \nabla\phi(c)\right),\\[3mm]
\tau c_t=\varepsilon \triangle c+ g(u,c),
\end{array}
\right.
\end{equation}
where $u(x,t)$ and $c(x,t)$ denote the cell density and the chemical
concentration, respectively. $D>0$ is the diffusion rate of cells
(bacteria) and $\varepsilon\geq0$ is the diffusion rate of chemical
substance. $\tau\geq0$ is a relaxation time scale and $\chi>0$
corresponds to attractive chemotaxis. Here $g(u,c)$ is a kinetic
function and $\phi(c)$ denoting a chemotactical sensitivity
function. With different choices of $g(u,c)$ and $\phi(c)$, many
results have been established
 in the literatures, cf.\cite{Biler09, Hortsmann03,Winkler10}.

As in \cite{Li111,Li10,Li11}, if we consider the model
(\ref{Keller-Segel1}) with $\tau=1$, $\phi(c)=\ln c$, $g(u,c)=-\alpha uc$,
the resulting model reads
\begin{equation}\label{Keller-Segel1*}
\left\{
\begin{array}{l}
u_t=\nabla\cdot \left(D\nabla u-\chi u \nabla\ln c\right),\\[3mm]
c_t=\varepsilon \triangle c-\alpha uc.
\end{array}
\right.
\end{equation}
The model (\ref{Equation@}) is derived from (\ref{Keller-Segel1*})
through the Hopf-Cole transformation
\begin{equation}\label{Hopf-Cole}
\begin{array}{l}
\displaystyle{\bf q^\varepsilon}=-{\frac{\nabla c}{c}}=-\nabla\ln
c,\ p^\varepsilon=u
\end{array}
\end{equation}
and scalings
\begin{equation}\label{scalings-1}
\begin{array}{l}
\displaystyle\tilde{t}=\alpha t,\ \
\tilde{x}=x\sqrt{\frac{\alpha}{\chi}},\ \ \tilde{\bf q}={\bf
q}\sqrt{\frac{\chi}{\alpha}},\ \ \widetilde{D}=\frac{D}{\chi},\ \
 \widetilde{\varepsilon}=\frac{\varepsilon}{\chi},
\end{array}
\end{equation}
where tilde has been dropped. When the diffusion of chemical
substance is so small that  it is negligible, i.e,
$\varepsilon\rightarrow0^+$, then the model (\ref{Keller-Segel1*})
becomes
\begin{equation}\label{Keller-Segel-1}
\left\{
\begin{array}{l}
u_t=\nabla\cdot \left(D\nabla u-\chi u\nabla \ln c\right),\\[3mm]
c_t=-\alpha uc.
\end{array}
\right.
\end{equation}
A version of system (\ref{Keller-Segel-1}) was proposed by Othmer
and Stevens in \cite{Othmer97} to describe the chemotactic movement
of particles where the chemicals are non-diffusible. The models
developed in \cite{Othmer97} have been studied in depth by Levine
and Sleeman in \cite{Levine97}. They gave some heuristic
understanding of some of these phenomena and investigated the
properties of solutions of a system of chemotaxis equation arising
in the theory of reinforced random walks. Y. Yang, H. Chen and W.A.
Liu in \cite{Yang01} studied the global existence and blow-up in a
finite-time of solutions for the case considered in \cite{Levine97},
respectively.  For the other results on (\ref{Keller-Segel-1}),
please refer to \cite{Hillen04, Levine01,Yang05} and references
therein.

Similar to the derivation of system (\ref{Equation@}), the system
(\ref {Keller-Segel-1}) can be converted into a system of
conservation laws as follows:
\begin{equation}\label{Equation}
\left\{
\begin{array}{l}
p_t-\nabla\cdot\left(p\bf{q}\right)
=D\triangle p,\\[3mm]
{\bf{q}}_t-\nabla p=0,
\end{array}
\right.
\end{equation}
with initial data
\begin{equation}\label{Equation-initial}
\left(p, {\bf q}\right)(x,0)=\left(p_0(x),{\bf
q}_0(x)\right)\rightarrow(p_\infty, 0)\ \ {\rm as}\ \ |x|\rightarrow
\infty.
\end{equation}
System (\ref{Equation}) has been studied by several authors. For one
dimension, the global well-posedness of smooth solution was obtained
in \cite{Zhang07, Guo09} with small initial data and large initial
data, respectively. For high dimensions, the global well-posedness
of smooth solution to (\ref{Equation}) was investigated in
\cite{Li111, Li112} for Cauchy problem and initial-boundary value
problem, respectively, where the initial data is required to be
small at least in $H^2$ norm. For other related results, such as
nonlinear stability of waves in one dimension and so on, please
refer to \cite{Li09,Li10,Li11,Zhang12,Zhang121} and references
therein.


Formally, the system (\ref{Equation@}) becomes (\ref{Equation}) when
we take $\varepsilon=0$. In fact, the investigation of the problem
of the zero viscosity limit is one of the challenging topics in
fluid dynamics and has been much more extensively investigated for
many other models, cf.
\cite{Chen08,Frid99,Jiang09,Ruan11,YGWang,Xin99}. However, to our
knowledge, there are few results on the system (\ref{Equation@}) in
this direction, cf. \cite{Peng12}. Our aim here is to prove
accurately that the solutions of (\ref {Equation}) converge to the
solutions of (\ref {Equation@}) as the chemical diffusion
$\varepsilon$ goes to zero.

To do this, we first establish the global well-posedness of
classical solutions to the Cauchy problem for the diffusive model
(\ref{Equation@})
 with smooth initial
data which is of small $L^2$ norm. Some {\it a priori} estimates
independent of $t$ and $\varepsilon$ are also obtained. Then, based on
these estimates, we get the global
existence of classical solutions to the Cauchy problem for the
non-diffusive model (\ref{Equation}) after passing to the limits
$\varepsilon\rightarrow0$. Finally, we derive the convergence rate
of the diffusive model toward the non-diffusive model.

Before stating the main results, we explain some notations.

\bigbreak\noindent\textbf{Notations:}\ \ $L^p=L^p({\Omega})$ $(1\leq
p\leq\infty)$ denotes usual Lebesgue space with the norm
$$
\begin{array}{c}
\bigbreak\displaystyle \|f\|_{L^p({\Omega})}=\left(\int_{\Omega}
|f(x)|^pdx\right)^{\frac{1}{p}}, \ \
1\leq p<\infty,\\
\bigbreak\displaystyle\int f dx=\int_{R^3} f dx.
\end{array}
$$
$H^l({\Omega}) \ (l\geq 0)$ denotes the usual $l$th-order Sobolev
space with the norm
$$
\|f\|_l=\left(\sum\limits^l_{j=0}\|\partial^j_xf\|^2\right)^{\frac{1}{2}},
$$
where $\Omega=R^3$, and $\|\cdot\|=\|\cdot\|_0=||\cdot\|_{L^2}$.

\bigskip

The main results in this paper can be stated as follows:
\begin{theorem}\label{theorem 1 }
 For given number $M>0$ (not necessarily
small), assume that initial data $(p_0,{\bf q}_0)$ satisfy
\begin{equation}\label{Existence1-in}
\begin{array}{l}
 \nabla\times{\bf{q}}_0=0,\ \ \ \|\nabla p_0 \|_{L^2}^2+\|\nabla\cdot {\bf{q}}_0 \|_{L^2}^2\le M,\ \ \
\left(p_0-p_{\infty}, {\bf{q}}_0\right)\in H^3,\ \ \ p_{\infty}>0.
 \end{array}
 \end{equation}
Then there exists a positive constant $\varepsilon_0$ depending on
$M$ and $p_{\infty}$ such that the Cauchy problem exists a unique
global solution in ${R}^3 \times (0, \infty)$,
 which satisfies
\begin{equation}\label{Existence1-c}
\left\{
\begin{array}{l}
\left(p^\varepsilon- p_{\infty}, {\bf q^\varepsilon}\right) \in
L^\infty([0,\infty),H^3),\\[3mm]
 \nabla p^\varepsilon\in L^2([0,\infty),H^3),\ \ \nabla{\bf q}^\varepsilon \in L^2([0,\infty),H^2),\\[3mm]
 \varepsilon\nabla{\bf q}^\varepsilon \in L^2([0,\infty),H^3)
 \end{array}
 \right.
\end{equation}
and
\begin{equation}\label{Existence3-c}
\begin{array}{ll}
\displaystyle\left\|{\bf q}^\varepsilon(t)\right\|^2_{H^3}+\left\|p^\varepsilon(t)-p_\infty\right\|^2_{H^3}\\[3mm]
+\displaystyle\int_0^t\left( \|\nabla p^\varepsilon(s)\|^2_{H^3}+
\varepsilon\|\nabla {\bf
q}^\varepsilon(s)\|^2_{H^3}\right)ds+\int_0^t
\|\nabla{\bf q}^\varepsilon(s)\|^2_{H^2}ds\leq C,\\[3mm]
\end{array}
\end{equation}
where $C$ is a positive constant independent of $\varepsilon$ and t,
provided that
\begin{equation}\label{Existence1-sm}
\begin{array}{l}
 \|{\bf{q}}_0 \|_{L^2}^2+ \|
p_0-p_\infty\|_{L^2}^2\leq\varepsilon_0.
 \end{array}
 \end{equation}
\end{theorem}

\bigskip

The last theorem is concerned with the convergence rate as well as
the global well-posedness of
(\ref{Equation})-(\ref{Equation-initial}).

\begin{theorem}\label{theorem 2 }
Suppose that $(p_0, {\bf q}_0)$  satisfies the assumptions in
Theorem \ref{theorem 1 }, then the Cauchy problem
(\ref{Equation})-(\ref{Equation-initial}) exists a unique global
solution in ${R}^3 \times (0, \infty)$, which satisfies
 \begin{equation}\label{Existence1-c}
\left\{
\begin{array}{l}
\left(p- p_{\infty}, {\bf q}\right) \in
L^\infty([0,\infty),H^3),\\[3mm]
 \nabla p\in L^2([0,\infty),H^3),\ \ \nabla{\bf q} \in
 L^2([0,\infty),H^2).
\end{array}
 \right.
\end{equation}
\end{theorem}
\noindent Furthermore,
\begin{eqnarray}\label{Th-con-cb}
\arraycolsep=1.5pt
\begin{array}[b]{rl}\displaystyle
\left\|\left(p^\varepsilon-p\right)(t)\right\|_{H^2}
+\left\|\left({\bf q}^\varepsilon-{\bf
q}\right)(t)\right\|_{H^2}\leq&\displaystyle
C\varepsilon^{\frac{1}{2}}.
\end{array}
\end{eqnarray}
 In particular,
\begin{eqnarray}\label{theorem 2.2 last}
\arraycolsep=1.5pt
\begin{array}[b]{rl}\displaystyle
\left\|\left(p^\varepsilon-p\right)(t)\right\|_{L^\infty}
+\left\|\left({\bf q}^\varepsilon-{\bf
q}\right)(t)\right\|_{L^\infty}\leq&\displaystyle
C\varepsilon^{\frac{1}{2}}.
\end{array}
\end{eqnarray} Here $C$ is a positive constant independent of $\varepsilon$ and $t$.
\\

\begin{remark} Notice that for the global existence of the solutions
to Cauchy problem (\ref{Equation})-(\ref{Equation-initial}), we only
assume that the $L^2$ norm of initial data is small. The initial
data can be of large oscillations with constant state at far field.
This is an improvement of {\cite{Li111}} where the initial data is
required to be small in $H^s(s>\frac{d}{2}+1)$ norm which implies
the oscillations are small.
\end{remark}

\begin{remark}
The power of $\varepsilon$ in (\ref{Th-con-cb}) could be improved to
1, which needs a slightly modification of the proof of Theorem
\ref{theorem 2 } with more regular initial data. But in this case,
it seems that the coefficient $C$ might depend on $t$.
\end{remark}
 The proofs of Theorems \ref{theorem 1 } and \ref{theorem 2 } are based on the
classical energy method. The key point for the proof of Theorem
\ref{theorem 1 } is to obtain some {\it a priori} estimates independent of $\varepsilon$ in which the $L^2-$bound of $\nabla\cdot {\bf
q}^\varepsilon$ plays a crucial role. For the proof of Theorem
\ref{theorem 2 }, some estimates of the order
$O\left(\varepsilon^{\frac{1}{2}}\right)$ are required, which needs
some delicate analysis.

The rest of the paper is organized as follows. In Section 2, we
study the global unique solvability on the Cauchy problem
(\ref{Equation@})-(\ref{Initial}). In Section 3, the zero-diffusion
limit as well as the global well-posedness of the solutions to
(\ref{Equation})-(\ref{Equation-initial}) is considered. We show
that the convergence rate in $L^\infty$-norm is of the order
$O(\varepsilon^{\frac{1}{2}})$, when diffusion parameter
$\varepsilon\rightarrow0^+$.

Throughout this paper, we denote a generic positive constant by $C$
which is independent of $\varepsilon$ and $t$.


\section{Proof of Theorem \ref{theorem 1 }}
\setcounter{equation}{0} In this section, we are concerned with the
global existence of large-oscillations solutions to the Cauchy
problem (\ref{Equation@})-(\ref{Initial}) when the initial data is
sufficiently close to a constant in $L^2$-norm. The global existence
 follows from a local existence theorem and some {\it a priori} estimates globally in time.

The local existence of the solutions could be done by using some
arguments similar to \cite{Li111}. We shall get some {\it a priori}
estimates globally in time which are also uniform for $\varepsilon$.

More precisely, for any
 given $T>0$ and $\varepsilon\geq0$, suppose $\left(p^\varepsilon(x,t), {\bf
q}^\varepsilon(x,t)\right)$ is a smooth solution to the Cauchy
problem (\ref{Equation@})-(\ref{Initial}) with regularities as in
Theorem \ref{theorem 1 }, we get the following key proposition.
 \begin{proposition}\label{prop 1}  For given number $M>0$ (not necessarily
small), assume that initial data $(p_0,{\bf q}_0)$ satisfy
$$
 \nabla\times{\bf{q}}_0=0,\ \ \ \|\nabla p_0 \|_{L^2}^2+\|\nabla\cdot {\bf{q}}_0 \|_{L^2}^2\le M,\ \ \
\left(p_0-p_{\infty}, {\bf{q}}_0\right)\in H^3,\ \ \ p_{\infty}>0.
$$ Then there exists some positive constant $\varepsilon_0$ depending on
$M$ and $p_\infty$ such that if $\left(p^\varepsilon(x,t), {\bf
 q}^\varepsilon(x,t)\right)$ is a smooth solution to the Cauchy
 problem of
 (\ref{Equation@})-(\ref{Initial}) in ${\mathbb{R}^3}\times(0,T]$, satisfying
\begin{equation}\label{3priori assumption-c}
\|{\bf q}^\varepsilon\|^2_{L^2}\leq 2\varepsilon_0,\ \ \|\nabla\cdot
{\bf{q}}^\varepsilon \|_{L^2}^2\leq 2M,
\end{equation}
then
\begin{equation}\label{3result}
\|{\bf q}^\varepsilon\|^2_{L^2}\leq \frac{3}{2}\varepsilon_0,\ \
\|\nabla\cdot {\bf{q}}^\varepsilon \|_{L^2}^2\leq \frac{3}{2}M,
\end{equation} provided that $$
\|{\bf{q}}_0 \|_{L^2}^2+ \|
p_0-p_\infty\|_{L^2}^2\leq\varepsilon_0.$$
\end{proposition}

\bigskip

 Without loss of generality, we suppose that
$D=1, p_{\infty}=1$. Letting $\widetilde{p}=p-1$, we obtain that
\begin{equation}\label{3Equation@}
\left\{
\begin{array}{l}
\widetilde{p}^\varepsilon_t-\nabla\cdot\left(\widetilde{p}^\varepsilon
\bf{q}^\varepsilon\right)-\nabla\cdot{\bf{q}}^\varepsilon
=\triangle \widetilde{p}^\varepsilon,\\[3mm]
{\bf{q}}^\varepsilon_t+\nabla\left(\varepsilon
|{\bf{q}}^\varepsilon|^2-\widetilde{p}^\varepsilon\right)=\varepsilon
\triangle {\bf{q}}^\varepsilon,
\end{array}
\right.
\end{equation}
with initial data
\begin{equation}\label{Initial 1@}
\left(\widetilde{p}^\varepsilon, {\bf
q^\varepsilon}\right)(x,0)=\left(p_0-1, {\bf
q}_0\right)\rightarrow(0, 0)\ \ {\rm as}\ \ |x|\rightarrow \infty.
\end{equation}
\bigskip

{\noindent\bf \underline{Proof of Proposition \ref{prop 1}}:}

\bigskip

The proof of Proposition \ref{prop 1} consists of the following
Lemmas \ref{le-2.1} and \ref{le-2.2}.
\begin{lemma}\label{le-2.1} [$L^2$ estimate]
Under the conditions of Theorem {\ref{theorem 1 }}, it holds that
\begin{eqnarray}\label{39Le-v-first11-3.5}
\arraycolsep=1.5pt
\begin{array}{lll}
\displaystyle\left\| \widetilde{p}^\varepsilon \right\|^2+\left\|
{\bf q}^\varepsilon\right\|^2+ \displaystyle\int^t_0\left\|
\nabla\widetilde{p}^\varepsilon\right\|^2
ds+\varepsilon\int^t_0\left\| \nabla{\bf q}^\varepsilon\right\|^2
ds\leq \frac{3}{2}\varepsilon_0,
\end{array}
\end{eqnarray} provided that $\varepsilon_0$ is small enough.
\end{lemma}
\pf Multiplying the first equation in (\ref{3Equation@}) by
$2\widetilde{p}^\varepsilon$ and the second by $2{\bf
q}^\varepsilon$, summing up them and then integrating over
$\mathbb{R}^3\times {[0,t]}$, one gets after integration by parts
that
\begin{eqnarray}\label{3Le-v-first1-c}
\arraycolsep=1.5pt
\begin{array}{lll}
&\displaystyle\left\| \widetilde{p}^\varepsilon \right\|^2+\left\|
{\bf q}^\varepsilon\right\|^2+
 2\int^t_0\left\| \nabla\widetilde{p}^\varepsilon\right\|^2 ds+2\varepsilon\int^t_0\left\| \nabla{\bf q}^\varepsilon\right\|^2 ds\\[3mm]
 =&\displaystyle\left\| \widetilde{p}^\varepsilon_0 \right\|^2+\left\| {\bf q}^\varepsilon_0\right\|^2
 +2\int^t_0\int\widetilde{p}^\varepsilon {\bf q}^\varepsilon\cdot(\nabla\widetilde{p}^\varepsilon) dxds
 +4\varepsilon\int^t_0\int({\bf{q}}^\varepsilon)^T\cdot(\nabla{\bf{q}}^\varepsilon\cdot{\bf q}^\varepsilon) dxds.
\end{array}
\end{eqnarray}
Next, we shall estimate the last two terms in the right-hand side.
By Cauchy inequality, H\"{o}lder inequality, Sobolev inequality and
Gagliardo-Nirenberg inequality, we obtain
\begin{eqnarray}\label{3Le-v-fi1-cc@}
\arraycolsep=1.5pt
\begin{array}[b]{lll}
\displaystyle2\int^t_0\int\widetilde{p}^\varepsilon {\bf
q}^\varepsilon\cdot(\nabla\widetilde{p}^\varepsilon) dxds
 &\displaystyle\leq  \int^t_0\left\| \nabla\widetilde{p}^\varepsilon\right\|^2 ds+ C\int^t_0\left\| \widetilde{p}^\varepsilon{\bf q}^\varepsilon\right\|^2 ds\\[3mm]
&\displaystyle\leq   \int^t_0\left\| \nabla\widetilde{p}^\varepsilon\right\|^2 ds+ C\int^t_0\left\| \widetilde{p}^\varepsilon\|_{L^6}^2\|{\bf q}^\varepsilon\right\|^2_{L^3}ds\\[3mm]
&\displaystyle\leq   \int^t_0\left\| \nabla\widetilde{p}^\varepsilon\right\|^2 ds+ C\int^t_0\left\| \nabla\widetilde{p}^\varepsilon\right\|^2\|\nabla{\bf q}^\varepsilon\|\|{\bf q}^\varepsilon\|ds\\[3mm]
\end{array}
\end{eqnarray}
and
\begin{eqnarray}\label{3Le-v-fi-acc}
\arraycolsep=1.5pt
\begin{array}[b]{lll}
\displaystyle4\varepsilon\int^t_0\int({\bf{q}}^\varepsilon)^T\cdot(\nabla{\bf{q}}^\varepsilon\cdot
{\bf q}^\varepsilon) dxds
&\displaystyle\leq 4\varepsilon\int^t_0\left\|{\bf{q}}^\varepsilon\right\|_{L^3}\left\|\nabla{\bf{q}}^\varepsilon\right\|_{L^2} \left\|{\bf q}^\varepsilon\right\|_{L^6} ds\\[3mm]
&\displaystyle\leq C\varepsilon\int^t_0\|\nabla{\bf
q}^\varepsilon\|^{\frac{1}{2}}\|{\bf
q}^\varepsilon\|^{\frac{1}{2}}\left\|\nabla{\bf{q}}^\varepsilon\right\|^2ds.
\end{array}
\end{eqnarray}
Since $\triangle{\bf q}^\varepsilon= \nabla(\nabla\cdot{\bf
q}^\varepsilon)-\nabla\times(\nabla\times{\bf q}^\varepsilon)$, we
obtain that
\begin{eqnarray*}
\arraycolsep=1.5pt
\begin{array}[b]{rl}
\nabla{\bf q}^\varepsilon=
-\nabla(-\triangle)^{-1}\nabla(\nabla\cdot{\bf
q}^\varepsilon)+\nabla(-\triangle)^{-1}\nabla\times(\nabla\times{\bf
q}^\varepsilon).
\end{array}
\end{eqnarray*}
The standard $L^2$ estimate shows that
\begin{eqnarray}\label{div-curl}
\arraycolsep=1.5pt
\begin{array}[b]{lll}
\left\|\nabla{\bf q}^\varepsilon\right\|_{L^2}\leq C
\left(\left\|\nabla\cdot{\bf
q}^\varepsilon\right\|_{L^2}+\left\|\nabla\times{\bf
q}^\varepsilon\right\|_{L^2}\right).
\end{array}
\end{eqnarray}
Moreover, by taking the curl for $(\ref{3Equation@})_{2}$, one has
\begin{eqnarray}\label{div-curl 1}
\arraycolsep=1.5pt
\begin{array}[b]{rl}&\displaystyle
 \frac{d}{dt}\left( \nabla\times{\bf q}^\varepsilon\right)=\varepsilon\triangle(\nabla\times{\bf q}^\varepsilon).
\end{array}
\end{eqnarray}
Initial data is given as
\begin{eqnarray}\label{div-curl 2}
\begin{array}[b]{rl}
(\nabla\times{\bf q}^\varepsilon)\mid_{t=0}=\nabla\times{\bf q}_0=0.
\end{array}
\end{eqnarray}
By solving the initial value problem (\ref{div-curl 1} ) and
(\ref{div-curl 2} ), one has
\begin{eqnarray}\label{cur}
\begin{array}[b]{rl}
\nabla\times{\bf q}^\varepsilon=0,
\end{array}
\end{eqnarray}
which implies
\begin{eqnarray}\label{div-cur3}
\arraycolsep=1.5pt
\begin{array}[b]{lll}
\triangle{\bf q}^\varepsilon= \nabla(\nabla\cdot{\bf
q}^\varepsilon),\ \ \ \ \left\|\nabla{\bf
q}^\varepsilon\right\|_{L^2}\leq C\left\|\nabla\cdot{\bf
q}^\varepsilon\right\|_{L^2}.
\end{array}
\end{eqnarray}
The combination of (\ref {3priori assumption-c}),
(\ref{3Le-v-fi1-cc@}), (\ref{3Le-v-fi-acc}) and (\ref {div-cur3})
yields
\begin{eqnarray}\label{3Le-v-first11-cc@}
\arraycolsep=1.5pt
\begin{array}[b]{lll}
\displaystyle2\int^t_0\int\widetilde{p}^\varepsilon {\bf
q}^\varepsilon \cdot(\nabla\widetilde{p}^\varepsilon) dxds \leq
\int^t_0\left\| \nabla\widetilde{p}^\varepsilon\right\|^2 ds+
C\sqrt{M\varepsilon_0}\int^t_0\left\|
\nabla\widetilde{p}^\varepsilon\right\|^2ds,
\end{array}
\end{eqnarray}
and
\begin{eqnarray}\label{3Le-v-first1-acc}
\arraycolsep=1.5pt
\begin{array}[b]{lll}
\displaystyle4\varepsilon\int^t_0\int({\bf{q}}^\varepsilon)^T\cdot(\nabla{\bf{q}}^\varepsilon\cdot
{\bf q}^\varepsilon) dxds \displaystyle\leq
C\varepsilon(M\varepsilon_0)^\frac{1}{4}\int^t_0\left\|\nabla{\bf{q}}^\varepsilon\right\|^2ds.
\end{array}
\end{eqnarray}
Substituting (\ref{3Le-v-first11-cc@}) and (\ref{3Le-v-first1-acc})
into (\ref{3Le-v-first1-c}) and setting
$\displaystyle\varepsilon_0\leq\frac{1}{16C^4M}$, one may arrive at
(\ref{39Le-v-first11-3.5}), where (\ref{Existence1-sm}) has been
used. This completes the proof of Lemma \ref{le-2.1}.
\endpf

\begin{lemma}\label{le-2.2} [First-order energy estimate]
Under the conditions of Theorem {\ref{theorem 1 }}, it holds that
\begin{eqnarray}\label{3Le-v-first11-3.9}
\arraycolsep=1.5pt
\begin{array}{lll}
\displaystyle\left\| \nabla\widetilde{p}^\varepsilon
\right\|^2+\left\| \nabla\cdot{\bf q}^\varepsilon\right\|^2+
\displaystyle\int^t_0\left\| \nabla\cdot{\bf
q}^\varepsilon\right\|^2
 \displaystyle+\int^t_0\left(\left\|
\widetilde{p}^\varepsilon_t\right\|^2ds
+\varepsilon\left\|\triangle{\bf
q}^\varepsilon\right\|^2\right)ds\leq \frac{3M}{2},
\end{array}
\end{eqnarray} provided that $\varepsilon_0$ is small enough.
\end{lemma}
\pf
 Notice that
\begin{eqnarray}\label{3Le-v-firstp2-cc@}
\arraycolsep=1.5pt
\begin{array}[b]{rl}\displaystyle
\nabla\cdot{\bf q}^\varepsilon_t=&\displaystyle
\triangle\widetilde{p}^\varepsilon-\varepsilon\triangle|{\bf q}^\varepsilon|^2+\varepsilon\triangle(\nabla\cdot{\bf q}^\varepsilon)\\[3mm]
=&\displaystyle\widetilde{p}^\varepsilon_t-\nabla\cdot\left(\widetilde{p}^\varepsilon{\bf
q}^\varepsilon\right)-\nabla\cdot{\bf q}^\varepsilon
-\varepsilon\triangle|{\bf
q}^\varepsilon|^2+\varepsilon\nabla\cdot(\triangle{\bf
q}^\varepsilon),
\end{array}
\end{eqnarray}
where we have used (\ref{div-cur3}).

Multiplying
(\ref{3Le-v-firstp2-cc@}) by $2\nabla\cdot{\bf q}^\varepsilon$ and
integrating the resulting equation over $\mathbb{R}^3$, one obtains
after integration by parts that
\begin{eqnarray}\label{31Le-v-firstp2-3.11}
\arraycolsep=1.5pt
\begin{array}[b]{rl}&\displaystyle
 \frac{d}{dt}\left\| \nabla\cdot{\bf q}^\varepsilon\right\|^2 +2\left\| \nabla\cdot{\bf q}^\varepsilon\right\|^2
 +2\varepsilon\left\|\triangle{\bf q}^\varepsilon\right\|^2\\[3mm]
 =&\displaystyle2\int\widetilde{p}^\varepsilon_t\nabla\cdot{\bf q}^\varepsilon dx-2\int\nabla\cdot\left(\widetilde{p}^\varepsilon{\bf q}^\varepsilon\right)\nabla\cdot{\bf q}^\varepsilon dx+2\varepsilon\int\nabla|{\bf q}^\varepsilon|^2(\triangle{\bf q}^\varepsilon) dx.
\end{array}
\end{eqnarray}
Next, multiplying the first equation in (\ref{3Equation@}) by
$2\widetilde{p}^\varepsilon_t$, integrating the resulting equality
over $\mathbb{R}^3$ and using integration by parts, one has
\begin{eqnarray}\label{3Le-v-first13-cc}
\arraycolsep=1.5pt
\begin{array}[b]{rl}\displaystyle
 \frac{d}{dt}\left\| \nabla\widetilde{p}^\varepsilon\right\|^2 +2\left\| \widetilde{p}^\varepsilon_t\right\|^2=
2\int\widetilde{p}^\varepsilon_t\nabla\cdot{\bf q}^\varepsilon
dx+2\int\nabla\cdot\left(\widetilde{p}^\varepsilon{\bf
q}^\varepsilon\right){\widetilde{p}^\varepsilon_t}dx.
\end{array}
\end{eqnarray}
The combination of (\ref{31Le-v-firstp2-3.11}) with
(\ref{3Le-v-first13-cc}) yields
\begin{eqnarray}\label{39Le-v-first1-3.13}
\arraycolsep=1.5pt
\begin{array}[b]{lll}
&\displaystyle \frac{d}{dt}\left(\left\| \nabla\cdot{\bf q}^\varepsilon\right\|^2+\left\| \nabla\widetilde{p}^\varepsilon\right\|^2\right) +2\left\| \nabla\cdot{\bf q}^\varepsilon\right\|^2+2\left\| \widetilde{p}^\varepsilon_t\right\|^2+2\varepsilon\left\|\triangle{\bf q}^\varepsilon\right\|^2\\[3mm]
=&\displaystyle4\int\widetilde{p}^\varepsilon_t\nabla\cdot{\bf q}^\varepsilon dx+2\int\nabla\cdot\left(\widetilde{p}^\varepsilon{\bf q}^\varepsilon\right){\widetilde{p}^\varepsilon_t}dx\\[3mm]
&\displaystyle-2\int\nabla\cdot\left(\widetilde{p}^\varepsilon{\bf q}^\varepsilon\right)\nabla\cdot{\bf q}^\varepsilon dx+2\varepsilon\int\nabla|{\bf q}^\varepsilon|^2\cdot(\triangle{\bf q}^\varepsilon) dx\\[3mm]
=&\displaystyle\sum\limits_{i=1}^4J_i. \\[3mm]
\end{array}
\end{eqnarray}
For $J_1$, using (\ref{3Le-v-firstp2-cc@}), and integration by
parts, we have
\begin{eqnarray}\label{39Le-v-firstp2-I1-cc}
\arraycolsep=1.5pt
\begin{array}[b]{lll}\displaystyle
J_1&\displaystyle=\displaystyle 4\frac{d}{dt}\int\nabla\cdot{\bf q}^\varepsilon\widetilde{p}^\varepsilon dx-4\int\nabla\cdot{\bf q}^\varepsilon_t\widetilde{p}^\varepsilon dx\\[3mm]
&\displaystyle=4\left\|
\nabla\widetilde{p}^\varepsilon\right\|^2+4\frac{d}{dt}\int\nabla\cdot{\bf
q}^\varepsilon\widetilde{p}^\varepsilon dx
-4\varepsilon\int\nabla|{\bf q}^\varepsilon|^2\cdot(\nabla\widetilde{p}^\varepsilon) dx-4\varepsilon\int\triangle{\bf q}^\varepsilon\cdot(\nabla\widetilde{p}^\varepsilon) dx\\[3mm]
&\displaystyle=4\left\|
\nabla\widetilde{p}^\varepsilon\right\|^2+\sum\limits_{i=1}^3J_1^i.
 \end{array}
\end{eqnarray}
By Cauchy inequality, H\"{o}lder inequality, Sobolev inequality,
Gagliardo-Nirenberg inequality and (\ref{3priori assumption-c}), we
estimate $J_1^2$-$J_{1}^3$ as follows:

\begin{eqnarray}\label{3Le-v-J_1^2}
\arraycolsep=1.5pt
\begin{array}[b]{lll}
J_1^2&\displaystyle=-8\varepsilon\int{\bf q}^\varepsilon\cdot(\nabla{\bf q}^\varepsilon)\cdot(\nabla\widetilde{p}^\varepsilon) dx\\[3mm]
&\displaystyle\leq C\varepsilon\left\|{\bf{q}}^\varepsilon\right\|_{L^3}\left\|\nabla{\bf q}^\varepsilon\right\|_{L^6}\left\| \nabla\widetilde{p}^\varepsilon\right\|_{L^2}\\[3mm]
&\displaystyle\leq
C\varepsilon\left(M\varepsilon_0\right)^{\frac{1}{4}}\left\|\triangle{\bf{q}}^\varepsilon\right\|\left\| \nabla\widetilde{p}^\varepsilon\right\|\\[3mm]
&\displaystyle\le
C\varepsilon\left(M\varepsilon_0\right)^{\frac{1}{4}}\left\|\triangle{\bf{q}}^\varepsilon\right\|^2
+C\left\| \nabla\widetilde{p}^\varepsilon\right\|^2
\end{array}
\end{eqnarray}
and
\begin{eqnarray}\label{3Le-v-J_1^3}
\arraycolsep=1.5pt
\begin{array}[b]{lll}
J_1^3 \displaystyle \le
\frac{1}{2}\varepsilon\left\|\triangle{\bf{q}}^\varepsilon\right\|^2+C\left\|
\nabla\widetilde{p}^\varepsilon\right\|^2.
\end{array}
\end{eqnarray}
On the other hand, Cauchy inequality gives
\begin{eqnarray}\label{39Le-v-firstp2-I2-c}
\arraycolsep=1.5pt
\begin{array}[b]{lll}\displaystyle
 J_2+J_3 \leq\displaystyle \frac{1}{2}\left\|\nabla\cdot{\bf q}^\varepsilon\right\|^2+\frac{1}{2}
 \left\|\widetilde{p}^\varepsilon_t\right\|^2+16\left\| \nabla\cdot\left(\widetilde{p}^\varepsilon{\bf
 q}^\varepsilon\right)\right\|^2.
 \end{array}
\end{eqnarray}
Similar to (\ref {3Le-v-J_1^2}), it is immediate to obtain
\begin{eqnarray}\label{3Le-v-2.25}
\arraycolsep=1.5pt
\begin{array}[b]{lll}
J_4\displaystyle\leq
C\varepsilon(M\varepsilon_0)^{\frac{1}{4}}\left\|\triangle{\bf{q}}^\varepsilon\right\|^2.
\end{array}
\end{eqnarray}
Finally, we estimate the last term on the right hand side of
(\ref{39Le-v-firstp2-I2-c})
\begin{eqnarray}\label{39Le-v-2.26}
\arraycolsep=1.5pt
\begin{array}[b]{lll}\displaystyle
\left\| \nabla\cdot\left(\widetilde{p}^\varepsilon{\bf
q}^\varepsilon\right)\right\|^2\leq C\left\|
\nabla\widetilde{p}^\varepsilon\cdot{\bf
q}^\varepsilon\right\|^2+C\left\|
{\widetilde{p}^\varepsilon}\nabla\cdot{\bf
q}^\varepsilon\right\|^2=\sum\limits_{i=5}^6J_i.
 \end{array}
\end{eqnarray}
For $J_5$ and $J_6$, using Sobolev inequality, H\"{o}lder
inequality, Gagliardo-Nirenberg inequality, Young inequality,
(\ref{3priori assumption-c}) and (\ref{div-cur3}), we obtain
\begin{eqnarray}\label{39Le-v-2.27}
\arraycolsep=1.5pt
\begin{array}[b]{lll}
J_5&\displaystyle\leq C\left\| \nabla\widetilde{p}^\varepsilon\right\|_{L^6}^2\left\|{\bf q}^\varepsilon\right\|_{L^3}^2\\[3mm]
&\displaystyle\leq C\left\|
\nabla^2\widetilde{p}^\varepsilon\right\|^2\left\|\nabla{\bf
q}^\varepsilon\right\|\left\|{\bf q}^\varepsilon\right\|\leq
C\sqrt{M\varepsilon_0}\left\|
\nabla^2\widetilde{p}^\varepsilon\right\|^2
 \end{array}
\end{eqnarray}
and
\begin{eqnarray}\label{39Le-v-2.28}
\arraycolsep=1.5pt
\begin{array}[b]{rl}\displaystyle
J_6\leq&C\displaystyle
\left(\left\|\widetilde{p}^\varepsilon\right\|_{L^4}^2
 +\left\| \nabla\widetilde{p}^\varepsilon\right\|_{L^4}^2\right)\left\|\nabla\cdot{\bf q}^\varepsilon\right\|^2
 \\[3mm]
  \displaystyle
\leq&\displaystyle
 C\displaystyle \left(\left\|\widetilde{p}^\varepsilon\right\|^{\frac{1}{2}} \left\|\nabla\widetilde{p}^\varepsilon\right\|^{\frac{3}{2}}
 +\left\|\nabla\widetilde{p}^\varepsilon\right\|^{\frac{1}{2}} \left\|\nabla^2\widetilde{p}^\varepsilon\right\|^{\frac{3}{2}}\right)
 \left\|\nabla\cdot{\bf q}^\varepsilon\right\|^2\\[3mm]
\leq&\displaystyle C\left\|\nabla\widetilde{p}^\varepsilon\right\|^2\left\|\nabla\cdot{\bf q}^\varepsilon\right\|^2
+C\varepsilon_0\left\|\nabla\cdot{\bf q}^\varepsilon\right\|^2\\[3mm]
&\displaystyle+
C\left\|\nabla\widetilde{p}^\varepsilon\right\|^2\left\|\nabla\cdot{\bf
q}^\varepsilon\right\|^8
+\varepsilon_1\left\|\nabla^2\widetilde{p}^\varepsilon\right\|^2\\[3mm]
 \leq&\displaystyle CM(1+M^3)\left\|\nabla\widetilde{p}^\varepsilon\right\|^2
+\varepsilon_1\left\|\nabla^2\widetilde{p}^\varepsilon\right\|^2+C\varepsilon_0\left\|\nabla\cdot{\bf
q}^\varepsilon\right\|^2.
 \end{array}
\end{eqnarray}
Substituting (\ref{39Le-v-2.27})-(\ref{39Le-v-2.28}) into
(\ref{39Le-v-2.26}) shows that
\begin{eqnarray}\label{39Le-v-firstp2-I6-c}
\arraycolsep=1.5pt
\begin{array}[b]{lll}\displaystyle
\displaystyle \left\| \nabla\cdot\left(\widetilde{p}^\varepsilon{\bf
q}^\varepsilon\right)\right\|^2 \leq \displaystyle
\left(C\sqrt{M\varepsilon_0}+\varepsilon_1\right)\left\|\nabla^2\widetilde{p}^\varepsilon\right\|^2
+CM(1+M^3)\left\|\nabla\widetilde{p}^\varepsilon\right\|^2+
C\varepsilon_0\left\|\nabla\cdot{\bf q}^\varepsilon\right\|^2.
\end{array}
\end{eqnarray}
Applying the standard $L^2$-estimate for $(\ref{3Equation@})_1$, one
has
\begin{eqnarray}\label{39Le-v-firstp2-I4-es}
\arraycolsep=1.5pt
\begin{array}[b]{lll}\displaystyle
\left\| \nabla^2\widetilde{p}^\varepsilon\right\|^2\leq
C\left(\left\|\nabla\cdot{\bf
q}^\varepsilon\right\|^2+\displaystyle\left\|\widetilde{p}^\varepsilon_t\right\|^2+\left\|
\nabla\cdot\left(\widetilde{p}^\varepsilon{\bf
q}^\varepsilon\right)\right\|^2\right),
 \end{array}
\end{eqnarray}
which together with (\ref{39Le-v-firstp2-I6-c}) gives
\begin{eqnarray}\label{39Le-v-firstp2-I7-c}
\arraycolsep=1.5pt
\begin{array}[b]{lll}\displaystyle
\left\| \nabla\cdot\left(\widetilde{p}^\varepsilon{\bf
q}^\varepsilon\right)\right\|^2
 \leq&\displaystyle \left(C\sqrt{M\varepsilon_0}+C\varepsilon_1\right)\left(\left\|\nabla\cdot{\bf q}^\varepsilon\right\|^2+\displaystyle\left\|\widetilde{p}^\varepsilon_t\right\|^2+\left\| \nabla\cdot\left(\widetilde{p}^\varepsilon{\bf q}^\varepsilon\right)\right\|^2\right)\\[3mm]
&\displaystyle+
CM(1+M^3)\left\|\nabla\widetilde{p}^\varepsilon\right\|^2+
C\varepsilon_0\left\|\nabla\cdot{\bf q}^\varepsilon\right\|^2.
 \end{array}
\end{eqnarray}
Taking $\varepsilon_1$ suitably small, one can deduce that
\begin{eqnarray}\label{2.32}
\arraycolsep=1.5pt
\begin{array}[b]{lll}\displaystyle
\left\| \nabla\cdot\left(\widetilde{p}^\varepsilon{\bf
q}^\varepsilon\right)\right\|^2
 \leq &\displaystyle\left(C\sqrt{M\varepsilon_0}+C\varepsilon_1\right)\left(\left\|\nabla\cdot{\bf q}^\varepsilon\right\|^2+\displaystyle\left\|\widetilde{p}^\varepsilon_t\right\|^2\right)\\[3mm]
&\displaystyle+
CM(1+M^3)\left\|\nabla\widetilde{p}^\varepsilon\right\|^2+
C\varepsilon_0\left\|\nabla\cdot{\bf q}^\varepsilon\right\|^2,
 \end{array}
 \end{eqnarray}
provided $\displaystyle\varepsilon_0\le\frac{(1-2C\varepsilon_1)^2}{4C^2M}$.\\

Substituting (\ref{39Le-v-firstp2-I1-cc})-(\ref{3Le-v-2.25}) and
(\ref{2.32}) into (\ref{39Le-v-first1-3.13}) and integrating the
resulting inequality over $[0,t)$, one gets after using
(\ref{Existence1-in}), (\ref{Existence1-sm}) and
(\ref{39Le-v-first11-3.5}) that
\begin{eqnarray*}
\arraycolsep=1.5pt
\begin{array}[b]{lll}
&\displaystyle \left\| \nabla\cdot{\bf q}^\varepsilon\right\|^2+
\left\| \nabla\widetilde{p}^\varepsilon\right\|^2 +2\int_0^t\left\|
\nabla\cdot{\bf q}^\varepsilon\right\|^2ds+2\int_0^t\left(\left\|
\widetilde{p}^\varepsilon_t\right\|^2
+\varepsilon\left\|\triangle{\bf q}^\varepsilon\right \|^2\right) ds\\[3mm]
\leq&\displaystyle
M+\left(C\sqrt{M\varepsilon_0}+C\varepsilon_1+C\varepsilon_0
+\frac{1}{2}\right)\int_0^t\left(\left\|\nabla\cdot{\bf
q}^\varepsilon\right\|^2
+\displaystyle\left\|\widetilde{p}^\varepsilon_t\right\|^2\right)ds\\[3mm]
&\displaystyle+
C\varepsilon_0\left(M+M^4+1\right)+\left(C(M\varepsilon_0)^{\frac{1}{4}}
+\frac{1}{2}\right)\varepsilon\int^t_0\left\|\triangle{\bf{q}}^\varepsilon\right\|^2ds\\[3mm]
&\displaystyle+ 4\int\nabla\cdot{\bf
q}^\varepsilon\widetilde{p}^\varepsilon dx-4\int\nabla\cdot{\bf
q}_0\widetilde{p}_0 dx,
\end{array}
\end{eqnarray*} which together with Cauchy inequality and (\ref{3priori
assumption-c}) deduces
\begin{eqnarray}\label{39Le-v-first1-I8-cc}
\arraycolsep=1.5pt
\begin{array}[b]{lll}
&\displaystyle \left\| \nabla\cdot{\bf q}^\varepsilon\right\|^2+
\left\| \nabla\widetilde{p}^\varepsilon\right\|^2 +2\int_0^t\left\|
\nabla\cdot{\bf q}^\varepsilon\right\|^2ds+2\int_0^t\left(\left\|
\widetilde{p}^\varepsilon_t\right\|^2
+\varepsilon\left\|\triangle{\bf q}^\varepsilon\right \|^2\right) ds\\[3mm]
\leq&\displaystyle
\frac{5M}{4}+\left(C\sqrt{M\varepsilon_0}+C\varepsilon_1+C\varepsilon_0
+\frac{1}{2}\right)\int_0^t\left(\left\|\nabla\cdot{\bf
q}^\varepsilon\right\|^2
+\displaystyle\left\|\widetilde{p}^\varepsilon_t\right\|^2\right)ds\\[3mm]
&\displaystyle+
C\varepsilon_0\left(M+M^4+1\right)+\left(C(M\varepsilon_0)^{\frac{1}{4}}
+\frac{1}{2}\right)\varepsilon\int^t_0\left\|\triangle{\bf{q}}^\varepsilon\right\|^2ds.
\end{array}
\end{eqnarray}
Next, choosing $\varepsilon_1$ suitably small and taking
$$ \varepsilon_0\leq\displaystyle\min\left\{\frac{1}{16C^4M},\frac{M}{4C(M+M^4+1)},\frac{1-2C\varepsilon_1}{2C(\sqrt{M}+1)},1\right\},$$
 one can get (\ref{3Le-v-first11-3.9}). This completes the proof of Lemma \ref{le-2.2}.
\endpf

The proof of Proposition \ref{prop 1} is complete. To prove Theorem
\ref{theorem 1 }, we need to get some high order estimates. Before
beginning, we give the following corollary.

\begin{corollary}\label{coro 2.1}Under the conditions of Theorem {\ref{theorem 1 }}, it holds that
\begin{eqnarray}\label{39Le-v-firstp2-3.26}
\arraycolsep=1.5pt
\begin{array}[b]{lll}\displaystyle
\int_0^t\left(\|
\nabla^2\widetilde{p}^\varepsilon\|^2+\|\widetilde{p}^\varepsilon\|_{L^\infty}^2\right)ds\leq
C,
 \end{array}
\end{eqnarray} provided that $\varepsilon_0$ is small enough.
\end{corollary}
\pf It follows from (\ref{39Le-v-firstp2-I4-es}),
(\ref{39Le-v-firstp2-I7-c}) and (\ref{3Le-v-first11-3.9}) that
\begin{eqnarray}\label{39Le-v-firstp2-3.25}
\arraycolsep=1.5pt
\begin{array}[b]{lll}\displaystyle
\int_0^t\left\| \nabla^2\widetilde{p}^\varepsilon\right\|^2ds
&\displaystyle\leq C\left(\int_0^t\left\|\nabla\cdot{\bf
q}^\varepsilon\right\|^2ds+\int_0^t\displaystyle\left\|\widetilde{p}^\varepsilon_t\right\|^2ds+\int_0^t\left\|
\nabla\cdot\left(\widetilde{p}^\varepsilon{\bf
q}^\varepsilon\right)\right\|^2ds\right)\\[3mm]
&\leq C.
 \end{array}
\end{eqnarray}
By Sobolev's embedding theorem, combining (\ref{39Le-v-first11-3.5})
and (\ref{39Le-v-firstp2-3.25}), we get (\ref{39Le-v-firstp2-3.26}).
\endpf

\begin{lemma}\label{le-2.3} [Second-order energy estimate]
Under the conditions of Theorem {\ref{theorem 1 }}, it holds that
\begin{eqnarray}\label{3Le-v-second2-c}
\begin{array}[b]{ll}
\displaystyle&\left\| \nabla^2\widetilde{p}^\varepsilon
\right\|^2+\left\| \nabla(\nabla\cdot{\bf q}^\varepsilon)\right\|^2+
\displaystyle\int^t_0\left\| \nabla(\nabla\cdot{\bf q}^\varepsilon)\right\|^2 ds\\[3mm]
&\displaystyle+\int^t_0\left(\left\|
\nabla\widetilde{p}^\varepsilon_t\right\|^2
+\varepsilon\left\|\nabla\cdot(\triangle{\bf
q}^\varepsilon)\right\|^2\right)ds\leq C,
\end{array}
\end{eqnarray}provided that $\varepsilon_0$ is small enough.
\end{lemma}

\pf Differentiating (\ref{3Le-v-firstp2-cc@}) yields
\begin{eqnarray}\label{2Le-v-frstp2-cc@}
\arraycolsep=1.5pt
\begin{array}[b]{rl}\displaystyle
\displaystyle\nabla(\nabla\cdot{\bf q}^\varepsilon_t)
=&\displaystyle
\nabla(\triangle\widetilde{p}^\varepsilon)-\varepsilon\nabla(\triangle|{\bf
q}^\varepsilon|^2)
+\varepsilon\nabla(\triangle(\nabla\cdot{\bf q}^\varepsilon))\\[3mm]
=&\displaystyle\nabla\widetilde{p}^\varepsilon_t-\nabla\left(\nabla\cdot\left(\widetilde{p}^\varepsilon{\bf
q}^\varepsilon\right)\right)-\nabla(\nabla\cdot{\bf q}^\varepsilon)
-\varepsilon\nabla(\triangle|{\bf
q}^\varepsilon|^2)+\varepsilon\nabla(\nabla\cdot(\triangle{\bf
q}^\varepsilon)).
\end{array}
\end{eqnarray}
Multiplying (\ref{2Le-v-frstp2-cc@}) by $2\nabla(\nabla\cdot{\bf
q}^\varepsilon)$ and integrating the resulting equation over
$\mathbb{R}^3$, one obtains after integration by parts that
\begin{eqnarray}\label{31Le-v-frstp2-3.29}
\arraycolsep=1.5pt
\begin{array}[b]{rl}&\displaystyle
 \frac{d}{dt}\left\| \nabla(\nabla\cdot{\bf q}^\varepsilon)\right\|^2 +2\left\|\nabla( \nabla\cdot{\bf q}^\varepsilon)\right\|^2
 +2\varepsilon\left\|\nabla\cdot(\triangle{\bf q}^\varepsilon)\right\|^2\\[3mm]
 =&\displaystyle2\int\nabla(\nabla\cdot{\bf q}^\varepsilon)\cdot(\nabla\widetilde{p}^\varepsilon_t) dx
-2\int\nabla\left(\nabla\cdot\left(\widetilde{p}^\varepsilon{\bf q}^\varepsilon\right)\right)\cdot\left(\nabla\left(\nabla\cdot{\bf q}^\varepsilon\right)\right) dx\\
&\displaystyle +2\varepsilon\int\triangle|{\bf
q}^\varepsilon|^2\nabla\cdot\left(\triangle{\bf
q}^\varepsilon\right) dx.
\end{array}
\end{eqnarray}
Next, applying $\nabla$ to $(\ref{3Equation@})_1$, multiplying it by
$2\nabla\widetilde{p}^\varepsilon_t$, integrating the resulting
equality over $\mathbb{R}^3$ and using integration by parts, one has
\begin{eqnarray}\label{Le-v-first23-3.30}
\arraycolsep=1.5pt
\begin{array}[b]{rl}\displaystyle
 \frac{d}{dt}\left\| \triangle\widetilde{p}^\varepsilon\right\|^2 +2\left\| \nabla\widetilde{p}^\varepsilon_t\right\|^2=
2\int\nabla(\nabla\cdot{\bf
q}^\varepsilon)\cdot(\nabla\widetilde{p}^\varepsilon_t)dx+2\int\nabla(\nabla\cdot\left(\widetilde{p}^\varepsilon{\bf
q}^\varepsilon\right))\cdot(\nabla{\widetilde{p}^\varepsilon_t})dx.
\end{array}
\end{eqnarray}
Putting (\ref{31Le-v-frstp2-3.29}) and (\ref{Le-v-first23-3.30})
together, we have
\begin{eqnarray}\label{2.41}
\arraycolsep=1.5pt
\begin{array}[b]{lll}
&\displaystyle \frac{d}{dt}\left(\left\| \nabla(\nabla\cdot{\bf q}^\varepsilon)\right\|^2
+\left\| \triangle\widetilde{p}^\varepsilon\right\|^2\right) +2\left\| \nabla(\nabla\cdot{\bf q}^\varepsilon)\right\|^2+2\left\| \nabla\widetilde{p}^\varepsilon_t\right\|^2+2\varepsilon\left\|\nabla(\triangle{\bf q}^\varepsilon)\right\|^2\\[3mm]
=&\displaystyle4\int\cdot\nabla(\nabla\cdot{\bf q}^\varepsilon)\cdot(\nabla\widetilde{p}^\varepsilon_t) dx+2\int\nabla(\nabla\cdot\left(\widetilde{p}^\varepsilon{\bf q}^\varepsilon\right))\cdot(\nabla{\widetilde{p}^\varepsilon_t})dx\\[3mm]
&\displaystyle-2\int\nabla(\nabla\cdot\left(\widetilde{p}^\varepsilon{\bf q}^\varepsilon\right))\cdot\left(\nabla(\nabla\cdot{\bf q}^\varepsilon)\right)dx+2\varepsilon\int\triangle|{\bf q}^\varepsilon|^2\nabla\cdot(\triangle{\bf q}^\varepsilon) dx\\[3mm]
=&\displaystyle\sum\limits_{i=7}^{10}J_i. \\[3mm]
\end{array}
\end{eqnarray}
For $J_7$, using (\ref{2Le-v-frstp2-cc@}), and integration by parts,
we have
\begin{eqnarray}\label{2.42}
\arraycolsep=1.5pt
\begin{array}[b]{rll}\displaystyle
J_7\displaystyle=&\displaystyle
4\frac{d}{dt}\int\nabla(\nabla\cdot{\bf
q}^\varepsilon)\cdot(\nabla\widetilde{p}^\varepsilon)dx
-4\int\nabla({\nabla\cdot{\bf q}^\varepsilon})_t\cdot(\nabla\widetilde{p}^\varepsilon) dx\\[3mm]
=&\displaystyle 4\left\| \triangle\widetilde{p}^\varepsilon\right\|^2
+4\frac{d}{dt}\int\nabla(\nabla\cdot{\bf q}^\varepsilon)\cdot(\nabla\widetilde{p}^\varepsilon)dx\\[3mm]
&\displaystyle-4\varepsilon\int\triangle|{\bf
q}^\varepsilon|^2\triangle\widetilde{p}^\varepsilon dx
+4\varepsilon\int\nabla\cdot(\triangle{\bf q}^\varepsilon)\triangle\widetilde{p}^\varepsilon dx\\[3mm]
=&\displaystyle 4\left\|\triangle\widetilde{p}^\varepsilon\right\|^2+\sum\limits_{i=1}^3J_7^i.
 \end{array}
\end{eqnarray}
For $J_7^2$, using Sobolev inequality, H\"{o}lder inequality,
Gagliardo-Nirenberg inequality, Young inequality, Proposition
\ref{prop 1} and (\ref{div-cur3}), we obtain
\begin{eqnarray}\label{3Le-v-J_7^2}
\arraycolsep=1.5pt
\begin{array}[b]{lll}
J_7^2&\displaystyle=-8\varepsilon\int\left(|\nabla{\bf
q}^\varepsilon|^2
+\triangle {\bf q}^\varepsilon\cdot{\bf q}^\varepsilon\right)\triangle\widetilde{p}^\varepsilon dx \\[3mm]
&\displaystyle\leq
\frac{\varepsilon}{2}\left\|\triangle\widetilde{p}^\varepsilon\right\|^2
+C\varepsilon\left(\left\|\nabla{\bf q}^\varepsilon\right\|_{L^4}^2
+\left\|{\bf q}^\varepsilon\triangle {\bf q}^\varepsilon\right\|^2\right)\\[3mm]
&\displaystyle\leq
\frac{\varepsilon}{2}\left\|\triangle\widetilde{p}^\varepsilon\right\|^2
+C\varepsilon\left(\left\|\nabla{\bf
q}^\varepsilon\right\|^{\frac{1}{2}}\left\|\nabla^2{\bf
q}^\varepsilon\right\|^{\frac{3}{2}}
+\left\|\triangle {\bf q}^\varepsilon\right\|_{L^6}^2\left\|{\bf q}^\varepsilon\right\|_{L^3}^2\right)\\[3mm]
&\displaystyle\leq
\frac{\varepsilon}{2}\left\|\triangle\widetilde{p}^\varepsilon\right\|^2
+C\varepsilon\left\|\triangle{\bf
q}^\varepsilon\right\|^2+C\varepsilon\left\|\nabla{\bf
q}^\varepsilon\right\|^2+C\varepsilon \left\|\nabla(\triangle{\bf
q}^\varepsilon)\right\|^2\left\|{\bf
q}^\varepsilon\right\|\left\|\nabla{\bf
q}^\varepsilon\right\|\\[3mm]
&\displaystyle\leq
\frac{\varepsilon}{2}\left\|\triangle\widetilde{p}^\varepsilon\right\|^2
+C\varepsilon\left\|\triangle{\bf
q}^\varepsilon\right\|^2+C\varepsilon\left\|\nabla\cdot{\bf
q}^\varepsilon\right\|^2+C\varepsilon\sqrt{M\varepsilon_0}\left\|\nabla\cdot(\triangle{\bf
q}^\varepsilon)\right\|^2,
\end{array}
\end{eqnarray}
where in the last inequality we have used the following  fact
\begin{eqnarray*}
\arraycolsep=1.5pt
\begin{array}[b]{lll}
\displaystyle\left\|\nabla(\triangle {\bf q}^\varepsilon)\right\|^2
&\displaystyle\le C\left(\left\|\nabla\cdot(\triangle {\bf
q}^\varepsilon)\right\|^2
+\left\|\nabla\times(\triangle {\bf q}^\varepsilon)\right\|^2\right)\\[3mm]
&\displaystyle= C\left(\left\|\nabla\cdot(\triangle {\bf
q}^\varepsilon)\right\|^2
+\left\|\triangle(\nabla\times {\bf q}^\varepsilon)\right\|^2\right)\\[3mm]
&\displaystyle=C\left\|\nabla\cdot(\triangle {\bf
q}^\varepsilon)\right\|^2
\end{array}
\end{eqnarray*}
due to (\ref{div-curl}) and (\ref{cur}). Cauchy inequality gives
\begin{eqnarray}\label{3Le-v-J_1^3}
\arraycolsep=1.5pt
\begin{array}[b]{lll}
J_7^3 \displaystyle \le
\frac{1}{2}\varepsilon\left\|\nabla\cdot(\triangle{\bf{q}}^\varepsilon)\right\|^2
+C\left\|\triangle\widetilde{p}^\varepsilon\right\|^2
\end{array}
\end{eqnarray}
and
\begin{eqnarray}\label{3.33}
\arraycolsep=1.5pt
\begin{array}[b]{lll}\displaystyle
 J_8+J_9 \leq\displaystyle \frac{1}{2}\left\|\nabla(\nabla\cdot{\bf q}^\varepsilon)\right\|^2
 +\frac{1}{2}\left\|\nabla\widetilde{p}^\varepsilon_t\right\|^2
 +16\left\| \nabla(\nabla\cdot\left(\widetilde{p}^\varepsilon{\bf q}^\varepsilon\right))\right\|^2.
 \end{array}
\end{eqnarray}
Similar to $J_7^2$, we estimate $J_{10}$ as follows:
\begin{eqnarray}\label{3Le-v-acc@}
\arraycolsep=1.5pt
\begin{array}[b]{lll}
J_{10} \leq&\displaystyle
\frac{\varepsilon}{2}\left\|\nabla\cdot(\triangle{\bf
q}^\varepsilon)\right\|^2+C\varepsilon\sqrt{M\varepsilon_0}\left\|\nabla\cdot(\triangle{\bf
q}^\varepsilon)\right\|^2\\[3mm]
 &+\displaystyle C\varepsilon\left\|\triangle{\bf
q}^\varepsilon\right\|^2+C\varepsilon\left\|\nabla\cdot{\bf
q}^\varepsilon\right\|^2.
\end{array}
\end{eqnarray}
Finally, we estimate the last term on the right hand side of
(\ref{3.33})
\begin{eqnarray}\label{3.35}
\arraycolsep=1.5pt
\begin{array}[b]{lll}
&\displaystyle
\left\| \nabla(\nabla\cdot\left(\widetilde{p}^\varepsilon{\bf q}^\varepsilon\right))\right\|^2\\[3mm]
\leq&\displaystyle C\left\| \nabla^2\widetilde{p}^\varepsilon\cdot{\bf q}^\varepsilon\right\|^2
+C\left\| {\nabla\widetilde{p}^\varepsilon}\nabla\cdot{\bf q}^\varepsilon\right\|^2+C\left\| {\nabla\widetilde{p}^\varepsilon}\cdot(\nabla{\bf q}^\varepsilon)\right\|^2+C\left\| {\widetilde{p}^\varepsilon}\nabla(\nabla\cdot{\bf q}^\varepsilon)\right\|^2\\[3mm]
=&\displaystyle\sum\limits_{i={11}}^{14}J_i.
 \end{array}
\end{eqnarray}
By  H\"{o}lder inequality, Sobolev inequality, Gagliardo-Nirenberg
inequality, Young inequality, Proposition \ref{prop 1} and
(\ref{div-cur3}), we estimate $J_{11}$-$J_{14}$ as follows:
\begin{eqnarray}\label{3.36}
\arraycolsep=1.5pt
\begin{array}[b]{lll}\displaystyle
J_{11}&\displaystyle \leq C\left\| \nabla^2\widetilde{p}^\varepsilon\right\|_{L^6}^2\left\|{\bf q}^\varepsilon\right\|_{L^3}^2\\[3mm]
&\displaystyle \leq C\left\|
\nabla^3\widetilde{p}^\varepsilon\right\|^2\left\|\nabla{\bf
q}^\varepsilon\right\|\left\|{\bf q}^\varepsilon\right\|\leq
C\sqrt{M\varepsilon_0}\left\|
\nabla^3\widetilde{p}^\varepsilon\right\|^2,
 \end{array}
\end{eqnarray}

\begin{eqnarray}\label{3.37}
\arraycolsep=1.5pt
\begin{array}[b]{rl}\displaystyle
 J_{12}+J_{13}\leq&C\displaystyle \left\|\nabla\widetilde{p}^\varepsilon\right\|_{L^{\infty}}^2
\left(\left\|\nabla{\bf
q}^\varepsilon\right\|^2+\left\|\nabla\cdot{\bf
q}^\varepsilon\right\|^2\right)
 \\[3mm]
\leq&C\displaystyle \left\|\nabla^2\widetilde{p}^\varepsilon\right\|
 \left\| \nabla^3\widetilde{p}^\varepsilon\right\|\left\|\nabla\cdot{\bf q}^\varepsilon\right\|^2
 \\[3mm]
  \displaystyle
\leq&\displaystyle C\left\|\nabla\cdot{\bf q}^\varepsilon\right\|^4
\left\|\nabla^2\widetilde{p}^\varepsilon\right\|^2+\varepsilon_1\left\|\nabla^3\widetilde{p}^\varepsilon\right\|^2
 \end{array}
\end{eqnarray}
and
\begin{eqnarray}\label{3.38}
\arraycolsep=1.5pt
\begin{array}[b]{lll}\displaystyle
J_{14}\leq C\left\|
\widetilde{p}^\varepsilon\right\|_{L^{\infty}}^2\left\|
\nabla(\nabla\cdot{\bf q}^\varepsilon)\right\|^2.
 \end{array}
\end{eqnarray}
Substituting (\ref{3.36})-(\ref{3.38}) into (\ref{3.35}) shows that
\begin{eqnarray}\label{3.39}
\arraycolsep=1.5pt
\begin{array}[b]{lll}\displaystyle
&\displaystyle \left\| \nabla(\nabla\cdot\left(\widetilde{p}^\varepsilon{\bf q}^\varepsilon\right))\right\|^2\\[3mm]
\leq& \displaystyle
\left(C\sqrt{M\varepsilon_0}+\varepsilon_1\right)\left\|\nabla^3\widetilde{p}^\varepsilon\right\|^2+C\left\|\nabla^2\widetilde{p}^\varepsilon\right\|^2+
C\left\|
\widetilde{p}^\varepsilon\right\|_{L^{\infty}}^2\left\|\nabla(\nabla\cdot{\bf
q}^\varepsilon)\right\|^2.
\end{array}
\end{eqnarray}
Applying the standard $H^1$-estimate for $(\ref{3Equation@})_1$
leads to
\begin{eqnarray}\label{3.40}
\arraycolsep=1.5pt
\begin{array}[b]{lll}\displaystyle
\left\| \nabla^2\widetilde{p}^\varepsilon\right\|_{H^1}^2\leq
C\left(\left\|\nabla\cdot{\bf q}^\varepsilon\right\|_{H^1}^2
+\displaystyle\left\|\widetilde{p}^\varepsilon_t\right\|_{H^1}^2
+\left\| \nabla\cdot\left(\widetilde{p}^\varepsilon{\bf
q}^\varepsilon\right)\right\|_{H^1}^2\right),
 \end{array}
\end{eqnarray}
which together with (\ref{3.39}) gives
\begin{eqnarray}\label{2.55}
\arraycolsep=1.5pt
\begin{array}[b]{lll}\displaystyle
\left\| \nabla(\nabla\cdot\left(\widetilde{p}^\varepsilon{\bf
q}^\varepsilon\right))\right\|^2 \leq
&\displaystyle\left(C\sqrt{M\varepsilon_0}+C\varepsilon_1\right)\left(\left\|\nabla\cdot{\bf
q}^\varepsilon\right\|_{H^1}^2
+\displaystyle\left\|\widetilde{p}^\varepsilon_t\right\|_{H^1}^2
+\left\| \nabla\cdot\left(\widetilde{p}^\varepsilon{\bf q}^\varepsilon\right)\right\|_{H^1}^2\right)\\[3mm]
 &\displaystyle+
C\left\|\nabla^2\widetilde{p}^\varepsilon\right\|^2+ C\left\|
\widetilde{p}^\varepsilon\right\|_{L^{\infty}}^2\left\|\nabla(\nabla\cdot{\bf
q}^\varepsilon)\right\|^2.
 \end{array}
\end{eqnarray}
Substituting (\ref{2.42})-(\ref{3Le-v-acc@}) and (\ref{2.55}) into
(\ref{2.41}) and integrating the resulting inequality over $[0,t)$,
one gets after using (\ref{Existence1-in}), Lemmas
\ref{le-2.1}-\ref{le-2.2} and Corollary \ref{coro 2.1}  that
\begin{eqnarray}\label{39Le-v-first3-I8-cc}
\arraycolsep=1.5pt
\begin{array}[b]{lll}
&\displaystyle \left\| \nabla(\nabla\cdot{\bf
q}^\varepsilon)\right\|^2
+\left\| \triangle\widetilde{p}^\varepsilon\right\|^2\\[3mm]
&\displaystyle +2\int_0^t\left\| \nabla(\nabla\cdot{\bf
q}^\varepsilon)\right\|^2ds +2\int_0^t\left\|
\nabla\widetilde{p}^\varepsilon_t\right\|^2ds
+2\varepsilon\int_0^t\left\|\nabla\cdot(\triangle{\bf q}^\varepsilon)\right\|^2ds\\[3mm]
\leq&\displaystyle C +\frac{1}{2}\left\|\nabla(\nabla\cdot{\bf
q}^\varepsilon)\right\|^2
+\int_0^t\left\| \widetilde{p}^\varepsilon\right\|_{L^{\infty}}^2\left\|\nabla(\nabla\cdot{\bf q}^\varepsilon)\right\|^2ds\\[3mm]
&\displaystyle + \varepsilon\int_0^t\left\|\nabla\cdot(\triangle{\bf
q}^\varepsilon)\right\|^2ds
+C\varepsilon\sqrt{M\varepsilon_0}\int_0^t\left\|\nabla\cdot(\triangle{\bf q}^\varepsilon)\right\|^2ds\\[3mm]
&\displaystyle +\left(C\sqrt{M\varepsilon_0}+C\varepsilon_1
+\frac{1}{2}\right)\int_0^t\left(\left\|\nabla(\nabla\cdot{\bf
q}^\varepsilon)\right\|^2
+\displaystyle\left\|\nabla\widetilde{p}^\varepsilon_t\right\|^2\right)ds
.
\end{array}
\end{eqnarray}
One obtains (\ref{3Le-v-second2-c}) by using Corollary \ref{coro
2.1} and Gronwall's inequality.  This completes the proof of Lemma
\ref{le-2.3}.
\endpf

\begin{corollary}\label{coro 2.2}
Under the conditions of Theorem {\ref{theorem 1 }}, it holds that
\begin{eqnarray}\label{39Le-v-firstp2-3.26+1}
\arraycolsep=1.5pt
\begin{array}[b]{lll}\displaystyle
\int_0^t\left(\|
\nabla^2\widetilde{p}^\varepsilon\|_{H^1}^2+\|\nabla\widetilde{p}^\varepsilon\|_{L^\infty}^2\right)ds\leq
C,
 \end{array}
\end{eqnarray} provided that $\varepsilon_0$ is small enough.
\end{corollary}
\pf
 It follows from
(\ref{3.40})-(\ref{2.55}), (\ref{3Le-v-first11-3.9}) and
(\ref{3Le-v-second2-c}) that
\begin{eqnarray}\label{39Le-v--2.55}
\arraycolsep=1.5pt
\begin{array}[b]{lll}&\displaystyle
\int_0^t\left\| \nabla^2\widetilde{p}^\varepsilon\right\|_{H^1}^2ds\\[3mm]
\leq &\displaystyle C\left(\int_0^t\left\|\nabla\cdot{\bf
q}^\varepsilon\right\|_{H^1}^2ds+\int_0^t\displaystyle\left\|\widetilde{p}^\varepsilon_t\right\|_{H^1}^2ds+\int_0^t\left\|
\nabla\cdot\left(\widetilde{p}^\varepsilon{\bf
q}^\varepsilon\right)\right\|_{H^1}^2ds\right)\leq C.
\end{array}
\end{eqnarray}
 By Sobolev's embedding theorem, we get
\begin{eqnarray*}
\arraycolsep=1.5pt
\begin{array}[b]{lll}\displaystyle
\int_0^t\left\|\nabla\widetilde{p}^\varepsilon\right\|_{L^{\infty}}^2ds\leq
C,
 \end{array}
\end{eqnarray*}
which together with (\ref{39Le-v--2.55}) leads to (\ref{39Le-v-firstp2-3.26+1}).
\endpf

\begin{lemma}\label{le-2.4} [Higher-order energy estimate]
Under the conditions of Theorem {\ref{theorem 1 }}, it holds that
\begin{eqnarray}\label{3Le-v-second3-c}
\arraycolsep=1.5pt
\begin{array}[b]{ll}
\displaystyle&\left\| \nabla^3\widetilde{p}^\varepsilon
\right\|^2+\left\| \nabla^2(\nabla\cdot{\bf
q}^\varepsilon)\right\|^2
\displaystyle+\int^t_0\left\| \nabla^2(\nabla\cdot{\bf q}^\varepsilon)\right\|^2 ds\\[3mm]
&\displaystyle+\int^t_0\left(\left\|
\nabla^2\widetilde{p}^\varepsilon_t\right\|^2
+\varepsilon\left\|\nabla^3(\nabla\cdot{\bf
q}^\varepsilon)\right\|^2\right)ds\leq C,
\end{array}
\end{eqnarray}provided that $\varepsilon_0$ is small enough.
\end{lemma}
\pf Applying $\triangle$ to (\ref{3Le-v-firstp2-cc@}), multiplying
it by $\triangle(\nabla\cdot{\bf q}^\varepsilon)$, taking
integrations in $x$ and using integration by parts, one gets
\begin{eqnarray}\label{39Le-v-first1-cc1}
\arraycolsep=1.5pt
\begin{array}[b]{lll}
&\displaystyle \frac{d}{dt}\left\| \triangle(\nabla\cdot{\bf
q}^\varepsilon)\right\|^2 \displaystyle+2\left\|
\triangle(\nabla\cdot{\bf q}^\varepsilon)\right\|^2
+2\varepsilon\left\|\triangle^2{\bf q}^\varepsilon\right\|^2\\[3mm]
=&\displaystyle2\int\triangle\widetilde{p}^\varepsilon_t\triangle(\nabla\cdot{\bf
q}^\varepsilon)dx
-2\int\triangle(\nabla\cdot\left(\widetilde{p}^\varepsilon{\bf q}^\varepsilon\right))\triangle(\nabla\cdot{\bf q}^\varepsilon)dx\\[3mm]
&\displaystyle+2\varepsilon\int\nabla(\triangle|{\bf
q}^\varepsilon|^2)\cdot\left(\triangle^2{\bf q}^\varepsilon\right) dx.
\end{array}
\end{eqnarray}
Similar to (\ref{39Le-v-first1-cc1}), one has
\begin{eqnarray}\label{39Le-v-first1-cc2}
\arraycolsep=1.5pt
\begin{array}[b]{lll}
\displaystyle \frac{d}{dt} \left\|
\nabla(\triangle\widetilde{p}^\varepsilon)\right\|^2 +2\left\|
\triangle\widetilde{p}^\varepsilon_t\right\|^2
=\displaystyle2\int\Delta\widetilde{p}^\varepsilon_t\triangle(\nabla\cdot{\bf
q}^\varepsilon)dx
+2\int\triangle(\nabla\cdot\left(\widetilde{p}^\varepsilon{\bf
q}^\varepsilon\right))\triangle{\widetilde{p}^\varepsilon_t}dx.
\end{array}
\end{eqnarray}
Putting the above two equalities together, we get
\begin{eqnarray}\label{39Le-v-first1-cc}
\arraycolsep=1.5pt
\begin{array}[b]{lll}
&\displaystyle \frac{d}{dt}\left(\left\| \triangle(\nabla\cdot{\bf
q}^\varepsilon)\right\|^2
+\left\| \nabla(\triangle\widetilde{p}^\varepsilon)\right\|^2\right)\\[3mm]
&\displaystyle+2\left\| \triangle(\nabla\cdot{\bf
q}^\varepsilon)\right\|^2
+2\left\| \triangle\widetilde{p}^\varepsilon_t\right\|^2+2\varepsilon\left\|\triangle^2{\bf q}^\varepsilon\right\|^2\\[3mm]
=&\displaystyle4\int\triangle\widetilde{p}^\varepsilon_t\triangle(\nabla\cdot{\bf
q}^\varepsilon)dx
+2\int\triangle(\nabla\cdot\left(\widetilde{p}^\varepsilon{\bf q}^\varepsilon\right))\triangle{\widetilde{p}^\varepsilon_t}dx\\[3mm]
&\displaystyle-2\int\triangle(\nabla\cdot\left(\widetilde{p}^\varepsilon{\bf
q}^\varepsilon\right))\triangle(\nabla\cdot{\bf q}^\varepsilon)dx
+2\varepsilon\int\nabla(\triangle|{\bf q}^\varepsilon|^2)\cdot\left(\triangle^2{\bf q}^\varepsilon\right) dx\\[3mm]
=&\displaystyle\sum\limits_{i={15}}^{18}J_i. \\[3mm]
\end{array}
\end{eqnarray}
In a manner similar to the  estimates of $J_{7}$-$J_{10}$,
$J_{15}$-$J_{18}$ can be bounded as follows:
\begin{eqnarray}\label{39Le-v-J-15}
\arraycolsep=1.5pt
\begin{array}[b]{rll}\displaystyle
J_{15}=&\displaystyle 4\frac{d}{dt}\int\triangle(\nabla\cdot{\bf
q}^\varepsilon)\triangle\widetilde{p}^\varepsilon dx
-4\int{\triangle(\nabla\cdot{\bf q}^\varepsilon)}_t\triangle\widetilde{p}^\varepsilon dx\\[3mm]
 =&\displaystyle 4\left\| \nabla(\triangle\widetilde{p}^\varepsilon\right)\|^2
+4\frac{d}{dt}\int\triangle(\nabla\cdot{\bf
q}^\varepsilon)\triangle\widetilde{p}^\varepsilon dx
\\[3mm]
&\displaystyle-4\varepsilon\int\triangle^2{\bf
q}^\varepsilon\cdot\left(\nabla(\triangle\widetilde{p}^\varepsilon)\right) dx
-4\varepsilon\int\nabla(\triangle|{\bf
q}^\varepsilon|^2)\cdot\left(\nabla(\triangle\widetilde{p}^\varepsilon)\right) dx
\\[3mm]
=&\displaystyle4\left\|\nabla(\triangle\widetilde{p}^\varepsilon)\right\|^2+\sum\limits_{i=1}^3J_{15}^i.
 \end{array}
\end{eqnarray}
Cauchy inequality gives
\begin{eqnarray}\label{39Le-v-J15_12}
\arraycolsep=1.5pt
\begin{array}[b]{lrl}\displaystyle
J_{15}^2\le &\displaystyle
\displaystyle\frac{1}{2}\varepsilon\left\|\triangle^2{\bf{q}}^\varepsilon\right\|^2
+C\left\|\nabla(\triangle\widetilde{p}^\varepsilon)\right\|^2.
 \end{array}
\end{eqnarray}
For $J_{15}^3$, we have
\begin{eqnarray}\label{3Le-v4-acc@}
\arraycolsep=1.5pt
\begin{array}[b]{lll}
J_{15}^3 &\displaystyle\leq
\frac{\varepsilon}{2}\left\|\nabla(\triangle\widetilde{p}^\varepsilon)\right\|^2
+C\varepsilon\left(\left\|\nabla{\bf
q}^\varepsilon\right\|_{L^\infty}^2\left\|\triangle {\bf
q}^\varepsilon\right\|^2
+\left\|\nabla(\triangle {\bf q}^\varepsilon)\right\|_{L^6}^2\left\|{\bf q}^\varepsilon\right\|_{L^3}^2\right)\\[3mm]
&\displaystyle\leq
\frac{\varepsilon}{2}\left\|\nabla(\triangle\widetilde{p}^\varepsilon)\right\|^2
+C\varepsilon\left\|\nabla^3{\bf
q}^\varepsilon\right\|\left\|\nabla^2{\bf
q}^\varepsilon\right\|\left\|\triangle{\bf q}^\varepsilon\right\|^2
+C\varepsilon\sqrt{M\varepsilon_0}\left\|\triangle^2{\bf
q}^\varepsilon \right\|^2\\[3mm]
&\displaystyle\leq
\frac{\varepsilon}{2}\left\|\nabla(\triangle\widetilde{p}^\varepsilon)\right\|^2
+C\varepsilon\left\|\triangle(\nabla\cdot{\bf
q}^\varepsilon)\right\|^2+C\varepsilon\left\|\triangle{\bf
q}^\varepsilon\right\|^2
+C\varepsilon\sqrt{M\varepsilon_0}\left\|\triangle^2{\bf
q}^\varepsilon \right\|^2.\\[3mm]
\end{array}
\end{eqnarray}
By Cauchy inequality, we get
\begin{eqnarray}\label{39Le-v-firstp3-I2-c}
\arraycolsep=1.5pt
\begin{array}[b]{lll}\displaystyle
 J_{16}+J_{17} \leq\displaystyle \frac{1}{2}\left\|\triangle(\nabla\cdot{\bf q}^\varepsilon)\right\|^2
 +\frac{1}{2}\left\|\triangle\widetilde{p}^\varepsilon_t\right\|^2+C\left\| \triangle(\nabla\cdot\left(\widetilde{p}^\varepsilon{\bf q}^\varepsilon\right))\right\|^2.
 \end{array}
\end{eqnarray}
Similar to $J_{15}^3$, $J_{18}$ can be estimated as follows:
\begin{eqnarray}\label{3Le-v4-acc@}
\arraycolsep=1.5pt
\begin{array}[b]{lll}
J_{18} &\displaystyle\leq
\frac{\varepsilon}{2}\left\|\triangle^2{\bf q}^\varepsilon\right\|^2
+C\varepsilon\left\|\triangle(\nabla\cdot{\bf
q}^\varepsilon)\right\|^2+C\varepsilon\left\|\triangle{\bf
q}^\varepsilon\right\|^2
+C\varepsilon\sqrt{M\varepsilon_0}\left\|\triangle^2{\bf
q}^\varepsilon \right\|^2.
\end{array}
\end{eqnarray}
Next, we estimate the last term on the right hand side of
(\ref{39Le-v-firstp3-I2-c})
\begin{eqnarray}\label{39Le-v-firstp5-I3-c}
\arraycolsep=1.5pt
\begin{array}[b]{lll}
&\displaystyle
\left\| \triangle(\nabla\cdot\left(\widetilde{p}^\varepsilon{\bf q}^\varepsilon\right))\right\|^2\\[3mm]
\leq&\displaystyle C\left\| \nabla(\triangle\widetilde{p}^\varepsilon)\cdot{\bf q}^\varepsilon\right\|^2+C\left\| \nabla^2\widetilde{p}^\varepsilon\cdot(\nabla{\bf q}^\varepsilon)\right\|^2+C\left\| {\triangle\widetilde{p}^\varepsilon}\nabla\cdot{\bf q}^\varepsilon\right\|^2+C\left\|\nabla({\nabla\cdot{\bf q}^\varepsilon})\cdot\left({\nabla\widetilde{p}^\varepsilon}\right)\right\|^2\\[3mm]
&\displaystyle+C\left\| {\nabla\widetilde{p}^\varepsilon}\cdot(\triangle{\bf q}^\varepsilon)\right\|^2+C\left\| {\widetilde{p}^\varepsilon}\triangle(\nabla\cdot{\bf q}^\varepsilon)\right\|^2\\[3mm]
=&\displaystyle\sum\limits_{i={19}}^{24}J_i.
 \end{array}
\end{eqnarray}
We estimate $J_{19}$-$J_{24}$ as follows:
\begin{eqnarray}\label{39Le-v-firstp3-I4-c}
\arraycolsep=1.5pt
\begin{array}[b]{lll}\displaystyle
J_{19}\leq C\left\|
\nabla(\triangle\widetilde{p}^\varepsilon)\right\|_{L^6}^2\left\|{\bf
q}^\varepsilon\right\|_{L^3}^2\leq C\left\|
\triangle^2\widetilde{p}^\varepsilon\right\|^2\left\|\nabla{\bf
q}^\varepsilon\right\|\left\|{\bf q}^\varepsilon\right\|\leq
C\sqrt{M\varepsilon_0}\left\|\triangle^2\widetilde{p}^\varepsilon\right\|^2,
 \end{array}
\end{eqnarray}

\begin{eqnarray}\label{39Le-v-firstp3-I5-c}
\arraycolsep=1.5pt
\begin{array}[b]{rl}\displaystyle
 J_{20}+J_{21}\leq&C\displaystyle \left\|\nabla^3\widetilde{p}^\varepsilon\right\|
\left\|
\nabla^4\widetilde{p}^\varepsilon\right\|\left\|\nabla\cdot{\bf
q}^\varepsilon\right\|^2\\[3mm]
\leq&\displaystyle C\left\|\nabla\cdot{\bf q}^\varepsilon\right\|^4
\left\|\nabla(\triangle\widetilde{p}^\varepsilon)\right\|^2+\varepsilon_1\left\|\triangle^2\widetilde{p}^\varepsilon\right\|^2\\[3mm]
\leq&\displaystyle
C\left\|\nabla(\triangle\widetilde{p}^\varepsilon)\right\|^2
+\varepsilon_1\left\|\triangle^2\widetilde{p}^\varepsilon\right\|^2,
 \end{array}
\end{eqnarray}
\begin{eqnarray}\label{39Le-v-firstp3-Ii5-c}
\arraycolsep=1.5pt
\begin{array}[b]{rl}\displaystyle
J_{22}+J_{23}\leq&C\displaystyle
\left\|\nabla\widetilde{p}^\varepsilon\right\|^2_{L^{\infty}}\left\|\nabla(\nabla\cdot{\bf
q}^\varepsilon)\right\|^2 \le
C\left\|\nabla\widetilde{p}^\varepsilon\right\|^2_{L^{\infty}}
\end{array}
\end{eqnarray}
and
\begin{eqnarray}\label{39Le-v-firstp3-J4-c}
\arraycolsep=1.5pt
\begin{array}[b]{lll}\displaystyle
J_{24}\leq C\left\|
\widetilde{p}^\varepsilon\right\|_{L^{\infty}}^2\left\|
\triangle(\nabla\cdot{\bf q}^\varepsilon)\right\|^2.
 \end{array}
\end{eqnarray}
Substituting (\ref{39Le-v-firstp3-I4-c})-(\ref{39Le-v-firstp3-J4-c})
into (\ref{39Le-v-firstp5-I3-c}) shows that
\begin{eqnarray}\label{39Le-v-firstp3-J6-c}
\arraycolsep=1.5pt
\begin{array}[b]{lll}\displaystyle
 \left\| \triangle(\nabla\cdot\left(\widetilde{p}^\varepsilon{\bf q}^\varepsilon\right))\right\|^2
\leq &\displaystyle
\left(C\sqrt{M\varepsilon_0}+\varepsilon_1\right)\left\|\triangle^2\widetilde{p}^\varepsilon\right\|^2
+C\left\|\nabla(\triangle\widetilde{p}^\varepsilon)\right\|^2\\[3mm]
&+\displaystyle C\left\|
\nabla\widetilde{p}^\varepsilon\right\|_{L^{\infty}}^2+C\left\|
\widetilde{p}^\varepsilon\right\|_{L^{\infty}}^2\left\|
\triangle(\nabla\cdot{\bf q}^\varepsilon)\right\|^2.
\end{array}
\end{eqnarray}
Applying the standard $H^2$-estimate for $(\ref{3Equation@})_1$
leads to
\begin{eqnarray}\label{39Le-v-firstp3-I4-es}
\arraycolsep=1.5pt
\begin{array}[b]{lll}\displaystyle
\left\| \triangle\widetilde{p}^\varepsilon\right\|_{H^2}^2\leq
C\left(\left\|\nabla\cdot{\bf
q}^\varepsilon\right\|_{H^2}^2+\displaystyle\left\|\widetilde{p}^\varepsilon_t\right\|_{H^2}^2+\left\|
\nabla\cdot\left(\widetilde{p}^\varepsilon{\bf
q}^\varepsilon\right)\right\|_{H^2}^2\right),
 \end{array}
\end{eqnarray}
which together with (\ref{39Le-v-firstp3-J6-c}) gives
\begin{eqnarray}\label{39Le-v-firstp3-I7-c}
\arraycolsep=1.5pt
\begin{array}[b]{lll}\displaystyle
C\left\| \triangle(\nabla\cdot\left(\widetilde{p}^\varepsilon{\bf
q}^\varepsilon\right))\right\|^2 \leq& \displaystyle
\left(C\sqrt{M\varepsilon_0}+C\varepsilon_1\right)\left(\left\|\nabla\cdot{\bf
q}^\varepsilon\right\|_{H^2}^2+\displaystyle\left\|\widetilde{p}^\varepsilon_t\right\|_{H^2}^2
+\left\| \nabla\cdot\left(\widetilde{p}^\varepsilon{\bf q}^\varepsilon\right)\right\|_{H^2}^2\right)\\[3mm]
&+\displaystyle
C\left\|\nabla(\triangle\widetilde{p}^\varepsilon)\right\|^2+
C\left\|
\widetilde{p}^\varepsilon\right\|_{L^{\infty}}^2\left\|\triangle(\nabla\cdot{\bf
q}^\varepsilon)\right\|^2+C\left\|
\nabla\widetilde{p}^\varepsilon\right\|_{L^{\infty}}^2.
 \end{array}
\end{eqnarray}
Substituting (\ref{39Le-v-J-15})-(\ref{39Le-v-firstp5-I3-c}) and
(\ref{39Le-v-firstp3-I7-c}) into (\ref{39Le-v-first1-cc}) and
integrating the resulting inequality over $[0,t)$, one obtains after
using (\ref{Existence1-in}), Lemmas \ref{le-2.1}-\ref{le-2.2}, Lemma \ref{le-2.3},
Corollary \ref{coro 2.1} and Corollary  \ref{coro 2.2} that
\begin{eqnarray}\label{39Le-v-first3-I8-cc}
\arraycolsep=1.5pt
\begin{array}[b]{lll}
&\displaystyle \left\| \triangle(\nabla\cdot{\bf q}^\varepsilon)\right\|^2+\left\|\nabla(\triangle\widetilde{p}^\varepsilon)\right\|^2\\[3mm]
&+\displaystyle2\int_0^t\left\| \triangle(\nabla\cdot{\bf
q}^\varepsilon)\right\|^2ds
+2\int_0^t\left(\left\| \triangle\widetilde{p}^\varepsilon_t\right\|^2+\varepsilon\left\|\triangle^2{\bf q}^\varepsilon\right\|^2\right)ds\\[3mm]
\leq &\displaystyle C+ \frac{1}{2}\left\|\triangle(\nabla\cdot{\bf q}^\varepsilon)\right\|^2+\left(C\sqrt{M\varepsilon_0}+C\varepsilon_1
+\frac{1}{2}\right)\int_0^t\left(\left\|\triangle(\nabla\cdot{\bf
q}^\varepsilon)\right\|^2
+\displaystyle\left\|\triangle\widetilde{p}^\varepsilon_t\right\|^2\right)ds\\[3mm]
&+\displaystyle C\int_0^t\left\|
\widetilde{p}^\varepsilon\right\|_{L_{\infty}}^2\left\|\triangle(\nabla\cdot{\bf
q}^\varepsilon)\right\|^2ds+\varepsilon\int_0^t\left\|
\triangle^2{\bf q}^\varepsilon\right\|^2ds
+C\varepsilon\sqrt{M\varepsilon_0}\int_0^t\left\|\triangle^2{\bf
q}^\varepsilon\right\|^2ds.
\end{array}
\end{eqnarray}
Using Corollary \ref{coro 2.1}, Gronwall's inequality and
(\ref{div-cur3}), one can immediately get (\ref{3Le-v-second3-c}).
This completes the proof of Lemma \ref{le-2.4}.
\endpf

\begin{corollary}
Under the conditions of Theorem {\ref{theorem 1 }}, it holds that
\begin{eqnarray}\label{39Le-v-firstp2-I9-es}
\arraycolsep=1.5pt
\begin{array}[b]{lll}\displaystyle
\int_0^t\| \nabla^2\widetilde{p}^\varepsilon\|_{H^2}^2ds\leq C,
 \end{array}
\end{eqnarray} provided that $\varepsilon_0$ is small enough.
\end{corollary}
\pf It follows from (\ref{3Le-v-first11-3.9}),
(\ref{3Le-v-second2-c}), (\ref{3Le-v-second3-c}) and
(\ref{39Le-v-firstp3-I4-es}) that
\begin{eqnarray*}
\arraycolsep=1.5pt
\begin{array}[b]{lll}\displaystyle
\int_0^t\left\| \nabla^2\widetilde{p}^\varepsilon\right\|_{H^2}^2ds
\leq \displaystyle C\left(\int_0^t\left\|\nabla\cdot{\bf
q}^\varepsilon\right\|_{H_2}^2ds+\int_0^t\displaystyle\left\|\widetilde{p}^\varepsilon_t\right\|_{H_2}^2ds+\int_0^t\left\|
\nabla\cdot\left(\widetilde{p}^\varepsilon{\bf
q}^\varepsilon\right)\right\|_{H_2}^2ds\right)\leq C.
\end{array}
\end{eqnarray*}

\endpf

{\noindent\bf\underline{Proof of Theorem \ref{theorem 1 }:}}\\

As a consequence of (\ref{39Le-v-first11-3.5}),
(\ref{3Le-v-first11-3.9}),  (\ref{3Le-v-second2-c}),
(\ref{3Le-v-second3-c}) and (\ref{39Le-v-firstp2-I9-es}), one obtains
\begin{equation*}
\arraycolsep=1.5pt
\begin{array}{ll}
\displaystyle&\left\|{\bf
q}^\varepsilon(t)\right\|^2_{H^3}\displaystyle+\left\|\widetilde{p}^\varepsilon(t)\right\|^2_{H^3}
\displaystyle+\int_0^t \left(\|\nabla
\widetilde{p}^\varepsilon(s)\|^2_{H^3}+\varepsilon\|\nabla {\bf
q}^\varepsilon(s)\|^2_{H^3}\right)ds\\[3mm]&\displaystyle+\int_0^t \|\nabla{\bf
q}^\varepsilon(s)\|^2_{H^2}ds\leq C,
\end{array}
\end{equation*}
where we have used (\ref{cur}) and the following fact:
\begin{eqnarray}\label{div-cur7}
\arraycolsep=1.5pt
\begin{array}[b]{lll}
\left\|\nabla{\bf q}^\varepsilon\right\|_{H^s}\leq C
\left(\left\|\nabla\cdot{\bf
q}^\varepsilon\right\|_{H^s}+\left\|\nabla\times{\bf
q}^\varepsilon\right\|_{H^s}\right), \ \ \ \text{for} \ \ s=1,2,3.
\end{array}
\end{eqnarray}

The uniqueness result can be proved by the method used in
\cite{Li112}, we thus omit the details for brevity. This completes
the proof of  Theorem \ref{theorem 1 }.

\section{Proof of Theorem \ref{theorem 2 }}
\setcounter{equation}{0}

From (\ref{Existence3-c}), one can easily get a unique smooth
solution $(p,{\bf q)}$ to (\ref{Equation})-(\ref{Equation-initial})
with regularities as in Theorem \ref{theorem 2 } after passing to
the limits $\varepsilon\rightarrow0$ (take subsequence if
necessary).

We will give the proof of the last part of Theorem \ref{theorem 2 },
i.e., the convergence rate as $\varepsilon\rightarrow0$. To do this,
it suffices to show the following Proposition \ref{Le-3.6}.

\begin{proposition}\label{Le-3.6}
Assume that the assumptions listed in Theorem {\ref{theorem 1 }}
are satisfied. Then there exists a positive constant $C$,
independent of $t$ and $\varepsilon$, such that
\begin{eqnarray}\label{Le-con1lemma-c}
\arraycolsep=1.5pt
\begin{array}[b]{rl}
\displaystyle&\left\|{\bf q}^\varepsilon-{\bf q}\right\|_{H^2}^2
+\left\|p^\varepsilon-p\right\|_{H^2}^2\\[3mm]
&\displaystyle+\varepsilon\int_0^t\left\|\nabla({\bf
q}^\varepsilon-{\bf q})\right\|_{H^2}^2
+\left\|\nabla(p^\varepsilon-p)\right\|_{H^2}^2ds\leq C\varepsilon.
\end{array}
\end{eqnarray}
\end{proposition}

\pf Denote $D=1,\ \ \widetilde{p}=p-1$. Then
(\ref{Equation})-(\ref{Equation-initial}) can be translated into the
following system
\begin{equation}\label{limit equation}
\left\{
\begin{array}{l}
\widetilde{p}_t-\nabla\cdot\left(\widetilde{p}
\bf{q}\right)-\nabla\cdot{\bf{q}}
=\triangle \widetilde{p},\\[3mm]
{\bf{q}}_t-\nabla\widetilde{p}=0
\end{array}
\right.
\end{equation}
with initial data
\begin{equation}\label{limit Initial}
\left(\widetilde{p}, {\bf
q}\right)(x,0)=\left(\widetilde{p}_0(x),{\bf
q}_0(x)\right)\rightarrow(0, 0)\ \ {\rm as}\ \ |x|\rightarrow
\infty.
\end{equation}
Setting
\begin{equation}\label{Le-conp1-c}
\mathbf{\psi}^\varepsilon={\bf q}^\varepsilon-{\bf q},\ \ \
\theta^\varepsilon=p^\varepsilon-p=\widetilde{p}^\varepsilon-\widetilde{p}.
\end{equation}
Then we deduce from (\ref{limit equation})-(\ref{limit Initial}) and
(\ref{3Equation@})-(\ref{Initial 1@}) that
$\left(\psi^\varepsilon,\theta^\varepsilon\right)(x,t)$ satisfy the
following Cauchy problem:
\begin{equation}\label{Le-conp-equation-c}
\left\{
\begin{array}{l}
\theta_t^\varepsilon-\nabla\cdot\left(\psi^\varepsilon
\widetilde{p}+{\bf q}\theta^\varepsilon\right)-\nabla\cdot\psi^\varepsilon=\triangle\theta^\varepsilon,\\[3mm]
\psi_t^\varepsilon+\nabla\left(\varepsilon \left({\bf
q}^\varepsilon\right)^2-\theta^\varepsilon\right) =\varepsilon
\triangle\psi^\varepsilon+\varepsilon\triangle{\bf q}
\end{array}
\right.
\end{equation}
with initial data
\begin{equation}\label{Le-conp-initial-c}
\left(\psi^\varepsilon, \theta^\varepsilon\right)(x,0)=\left(0,
0\right).
\end{equation}
\noindent\textbf{\underline{Step 1}.}

\bigskip
\noindent Multiplying the fist and second equations of
(\ref{Le-conp-equation-c}) by $2\theta^\varepsilon$ and
$2\psi^\varepsilon$ respectively, integrating the adding result with
respect $x$ and $t$ over $\mathbb{R}^3\times {[0,t]}$, we have
\begin{eqnarray}\label{Le-conp4-c}
\arraycolsep=1.5pt
\begin{array}[b]{rl}\displaystyle
 &\displaystyle\left\|\psi^\varepsilon\right\|^2+\left\|\theta^\varepsilon\right\|^2
+2\int_0^t\left(\varepsilon\left\|\nabla\psi^\varepsilon\right\|^2+\left\|\nabla\theta^\varepsilon\right\|^2\right)ds\\[3mm]
 =&\displaystyle
 -4\varepsilon\int_0^t\int{\bf q}^\varepsilon\cdot(\nabla{\bf q}^\varepsilon)\cdot\psi^\varepsilon dxds
+2\varepsilon\int_0^t\int\psi^\varepsilon\cdot
(\triangle{\bf q})dxds\\[3mm]
 &\displaystyle-2\int_0^t\int
\widetilde{p} \nabla\theta^\varepsilon\cdot\psi^\varepsilon dxds
 -2\int_0^t\int{\bf q}\cdot(\nabla\theta^\varepsilon)\theta^\varepsilon dxds \\[3mm]
 =&\displaystyle\sum\limits_{i=1}^{4}K_i.
\end{array}
\end{eqnarray}
By Cauchy inequality, H\"{o}lder inequality, Sobolev inequality,
Gagliardo-Nirenberg inequality  and Theorem \ref{theorem 1 }, we
obtain
\begin{eqnarray}\label{Le-conp-J8-c}
\arraycolsep=1.5pt
\begin{array}[b]{rl}
K_1 &\displaystyle\leq 2\varepsilon\int_0^t\left\|\nabla{\bf
q}^\varepsilon\right\|^2ds
+2\varepsilon\int_0^t\left\|{\bf q}^\varepsilon \psi^\varepsilon\right\|^2ds\\[3mm]
&\displaystyle\leq C\varepsilon+C\varepsilon\int_0^t\left\|{\bf q}^\varepsilon\right\|_{L^3}^2\left\|\psi^\varepsilon\right\|_{L^6}^2ds\\[3mm]
&\displaystyle\leq C
\varepsilon+C\varepsilon\sqrt{M\varepsilon_0}\int_0^t\left\|\nabla\psi^\varepsilon\right\|^2ds.
\end{array}
\end{eqnarray}
Integration by parts and Cauchy inequality implies
\begin{eqnarray}\label{Le-conp-J9-c}
\arraycolsep=1.5pt
\begin{array}[b]{rl}\displaystyle
K_2 \leq\displaystyle \varepsilon\int_0^t\left\|\nabla{\bf
q}\right\|^2ds +\varepsilon\int_0^t
\left\|\nabla\psi^\varepsilon\right\|^2ds
\end{array}
\end{eqnarray}
and
\begin{eqnarray}\label{Le-conp-J10-J12-c}
\arraycolsep=1.5pt
\begin{array}[b]{rl}\displaystyle
\displaystyle K_3+K_4 \leq&\displaystyle
2\int_0^t\left(\left\|\widetilde{p}\right\|_{L^\infty}^2+\left\|{\bf
q}\right\|_{L^\infty}^2\right)
\left(\left\|\psi^\varepsilon\right\|^2
+\left\|\theta^\varepsilon\right\|^2\right)ds
\\[3mm]&\displaystyle+\frac{3}{2}\int_0^t\left\|\nabla\theta^\varepsilon\right\|^2ds. \\[3mm]

\end{array}
\end{eqnarray}
Substituting (\ref{Le-conp-J8-c})-(\ref{Le-conp-J10-J12-c}) into
(\ref{Le-conp4-c}), we get
\begin{eqnarray}\label{Le-conp8-4.10}
\arraycolsep=1.5pt
\begin{array}[b]{rl}\displaystyle
\displaystyle&\left\|\psi^\varepsilon\right\|^2+\left\|\theta^\varepsilon\right\|^2
+\displaystyle\int_0^t\left(\varepsilon\left\|\nabla\psi^\varepsilon\right\|^2+\left\|\nabla\theta^\varepsilon\right\|^2\right)ds\\[3mm]
 &\leq\displaystyle C \varepsilon+C\int_0^t\left(\left\|\widetilde{p}\right\|_{L^\infty}^2+\left\|{\bf q}\right\|_{L^\infty}^2\right)
\left(\left\|\psi^\varepsilon\right\|^2+\left\|\theta^\varepsilon\right\|^2\right)ds.
\end{array}
\end{eqnarray}

\noindent\textbf{\underline{Step 2}.}

\bigskip
\noindent Multiplying the fist and second equations of
(\ref{Le-conp-equation-c}) by $-2\triangle\theta^\varepsilon$ and
$-2\triangle\psi^\varepsilon$ respectively, integrating the adding
result with respect $x$ and $t$ over $\mathbb{R}^3\times {[0,t]}$,
we have
\begin{eqnarray}\label{Le-conp9-3.71}
\arraycolsep=1.5pt
\begin{array}[b]{rl}\displaystyle
 &\displaystyle\left\|\nabla\psi^\varepsilon\right\|^2+\left\|\nabla\theta^\varepsilon\right\|^2
+2\int_0^t\left(\varepsilon\left\|\triangle\psi^\varepsilon\right\|^2
+\left\|\triangle\theta^\varepsilon\right\|^2\right)ds\\[3mm]
 =&\displaystyle 4\varepsilon\int_0^t\int{\bf q}^\varepsilon\cdot (\nabla{\bf q}^\varepsilon)\cdot\triangle\psi^\varepsilon dxds
-2\varepsilon\int_0^t\int\triangle\psi^\varepsilon\cdot
(\triangle{\bf q})dxds\\[3mm]
&\displaystyle-2\int_0^t\int\nabla\cdot\left(\psi^\varepsilon
\widetilde{p}
+{\bf q} \theta^\varepsilon\right) \triangle\theta^\varepsilon dxds\\[3mm]
 =&\displaystyle\sum\limits_{i=5}^{7}K_i.
\end{array}
\end{eqnarray}
Cauchy inequality leads to
\begin{eqnarray}\label{Le-conp-72-c}
\arraycolsep=1.5pt
\begin{array}[b]{rl}\displaystyle
K_{5} \leq&\displaystyle C\varepsilon\int_0^t\left\|{\bf
q}^\varepsilon\right\|^2_{L^\infty}\left\|\nabla{\bf
q}^\varepsilon\right\|^2ds
+\frac{\varepsilon}{2}\int_0^t\left\|\triangle\psi^\varepsilon\right\|^2ds\\[3mm]

\end{array}
\end{eqnarray}
and
\begin{eqnarray}\label{Le-conp-J15-c}
\arraycolsep=1.5pt
\begin{array}[b]{rl}\displaystyle
K_{6} \leq\displaystyle \varepsilon\int_0^t\left\|\triangle{\bf
q}^\varepsilon\right\|^2ds
+\varepsilon\int_0^t\left\|\triangle\psi^\varepsilon\right\|^2ds.
\end{array}
\end{eqnarray}
Straightforward calculations show that:
\begin{eqnarray}\label{Le-conp-J17-c}
\arraycolsep=1.5pt
\begin{array}[b]{rl}\displaystyle
K_{7}=&\displaystyle-2\int_0^t\int\nabla\cdot\psi^\varepsilon
\widetilde{p}\triangle\theta^\varepsilon dxds
 -2\int_0^t\psi^\varepsilon\cdot(\nabla \widetilde{p}^\varepsilon)\triangle\theta^\varepsilon dxds\\[3mm]
  &\displaystyle-2\int_0^t\int \nabla\cdot{\bf q} \theta^\varepsilon\triangle\theta^\varepsilon dxds
 -2\int_0^t\int{\bf q} \cdot(\nabla\theta^\varepsilon)\triangle\theta^\varepsilon dxds\\[3mm]
 =&\displaystyle\sum\limits_{i=1}^4K_{7}^i.
\end{array}
\end{eqnarray}
Here $K_{7}^1-K_{7}^4$ are estimated as follows:
\begin{eqnarray}\label{Le-conp-J171-c}
\arraycolsep=1.5pt
\begin{array}[b]{rl}\displaystyle
K_{7}^1 \leq&\displaystyle
\frac{1}{4}\int_0^t\left\|\triangle\theta^\varepsilon\right\|^2ds
+C\int_0^t\left\|\widetilde{p}\right\|_{L^\infty}^2\left\|\nabla\psi^\varepsilon\right\|^2ds,
\end{array}
\end{eqnarray}

\begin{eqnarray}\label{Le-conp-J172-c}
\arraycolsep=1.5pt
\begin{array}[b]{rl}\displaystyle
K_{7}^2 \leq&\displaystyle
\frac{1}{4}\int_0^t\left\|\triangle\theta^\varepsilon\right\|^2ds
+C\int_0^t\left\|\nabla{\bf q}^\varepsilon\right\|_{L^\infty}^2\left\|\psi^\varepsilon\right\|^2ds\\[3mm]
\end{array}
\end{eqnarray}
and
\begin{eqnarray}\label{Le-conp-J173-4-c}
\arraycolsep=1.5pt
\begin{array}[b]{rl}\displaystyle
K_{7}^3+K_{7}^4 \leq&\displaystyle
\frac{1}{4}\int_0^t\left\|\triangle\theta^\varepsilon\right\|^2ds
+C\int_0^t\left\|\nabla\cdot{\bf q}\right\|_{L^\infty}^2\left\|\theta^\varepsilon\right\|^2ds\\[3mm]
&\displaystyle+C\int_0^t\left\|{\bf q}\right\|_{L^\infty}^2\left\|\nabla\theta^\varepsilon\right\|^2ds.\\[3mm]
\end{array}
\end{eqnarray}
Substituting (\ref{Le-conp-72-c})-(\ref{Le-conp-J173-4-c}) into
(\ref{Le-conp9-3.71}), we get
\begin{eqnarray}\label{3Le-conpr-cc}
\arraycolsep=1.5pt
\begin{array}[b]{rl}\displaystyle
 &\displaystyle\left\|\nabla\psi^\varepsilon\right\|^2+\left\|\nabla\theta^\varepsilon\right\|^2
+\int_0^t\left(\varepsilon\left\|\triangle\psi^\varepsilon\right\|^2
+\left\|\triangle\theta^\varepsilon\right\|^2\right)ds\\[3mm]
 \leq&\displaystyle C \varepsilon+C\int_0^t\left(\left\|\widetilde{p}\right\|^2_{W^{1,\infty}}+\left\|{\bf q}\right\|^2_{W^{1,\infty}}\right)
\left(\left\|\psi^\varepsilon\right\|_{H^1}^2+\left\|\theta^\varepsilon\right\|_{H^1}^2\right)ds.
\end{array}
\end{eqnarray}

\noindent\textbf{\underline{Step 3}.}

\bigskip
\noindent Differentiating (\ref{Le-conp-equation-c}) yields
\begin{equation}\label{dLe-conp-equation-c}
\left\{
\begin{array}{l}
\nabla\theta_t^\varepsilon-\nabla\left(\nabla\cdot\left(\psi^\varepsilon
\widetilde{p}+{\bf q}\theta^\varepsilon\right)\right)-\triangle\psi^\varepsilon=\nabla(\triangle\theta^\varepsilon),\\[3mm]
\nabla\cdot\psi_t^\varepsilon+\triangle\left(\varepsilon \left({\bf
q}^\varepsilon\right)^2-\theta^\varepsilon\right) =\varepsilon
\nabla\cdot(\triangle\psi^\varepsilon)+\varepsilon\nabla\cdot(\triangle{\bf
q}).
\end{array}
\right.
\end{equation}
Multiplying the fist and second equations of
(\ref{dLe-conp-equation-c}) by
$-2\nabla(\triangle\theta^\varepsilon)$ and
$-2\nabla\cdot(\triangle\psi^\varepsilon)$ respectively, integrating
the adding result with respect $x$ and $t$ over $\mathbb{R}^3\times
{[0,t]}$, we have
\begin{eqnarray}\label{Le-conp9-c}
\arraycolsep=1.5pt
\begin{array}[b]{rl}\displaystyle
 &\displaystyle\left\|\triangle\psi^\varepsilon\right\|^2+\left\|\triangle\theta^\varepsilon\right\|^2
+2\int_0^t\left(\varepsilon\left\|\nabla\cdot(\triangle\psi^\varepsilon)\right\|^2
+\left\|\nabla(\triangle\theta^\varepsilon)\right\|^2\right)ds\\[3mm]
 =&\displaystyle 2\varepsilon\int_0^t\int\nabla\cdot(\triangle\psi^\varepsilon)\triangle|{\bf q}^\varepsilon|^2 dxds
-2\varepsilon\int_0^t\int\nabla\cdot(\triangle\psi^\varepsilon)
\nabla\cdot(\triangle{\bf q})dxds\\[3mm]
&-\displaystyle2\int_0^t\int\nabla\left(\nabla\cdot\left(\psi^\varepsilon
\widetilde{p}
+{\bf q} \theta^\varepsilon\right) \right)\cdot\left(\nabla\left(\triangle\theta^\varepsilon\right)\right) dxds\\[3mm]
 =&\displaystyle\sum\limits_{i=8}^{10}K_i.
\end{array}
\end{eqnarray}
Then, it follows from Cauchy inequality that
\begin{eqnarray}\label{Le-conp-J81-c}
\arraycolsep=1.5pt
\begin{array}[b]{rl}\displaystyle
K_{8}=&\displaystyle
4\varepsilon\int_0^t\int\nabla\cdot(\triangle\psi^\varepsilon)
\left|\nabla{\bf q}^\varepsilon\right|^2 dxds
+4\varepsilon\int_0^t\int\nabla\cdot(\triangle\psi^\varepsilon) {\bf q}^\varepsilon\cdot(\triangle{\bf q}^\varepsilon) dxds\\[3mm]
\leq&\displaystyle C\varepsilon \int_0^t\left(\left\|\nabla{\bf
q}^\varepsilon\right\|_{L^4}^2 +\displaystyle \left\|{\bf
q}\triangle{\bf q}^\varepsilon\right\|^2\right)ds
+\frac{\varepsilon}{2}\int_0^t \left\|\nabla\cdot(\triangle\psi^\varepsilon)\right\|^2ds\\[3mm]
\leq&\displaystyle C\varepsilon\int_0^t\left(\left\|{\bf
q}^\varepsilon\right\|^2_{L^\infty}+\left\|\nabla{\bf
q}^\varepsilon\right\|_{L^\infty}^2\right)\left(\left\|\nabla{\bf
q}^\varepsilon\right\|^2+
\left\|\triangle{\bf q}^\varepsilon\right\|^2\right)ds\\[3mm]
&+\displaystyle\frac{\varepsilon}{2}\int_0^t
\left\|\nabla\cdot(\triangle\psi^\varepsilon)\right\|^2ds\\[3mm]
\leq&\displaystyle C\varepsilon+\frac{\varepsilon}{2}\int_0^t
\left\|\nabla\cdot(\triangle\psi^\varepsilon)\right\|^2ds,

\end{array}
\end{eqnarray}
\begin{eqnarray}\label{2Le-conp-J15-c}
\arraycolsep=1.5pt
\begin{array}[b]{rl}\displaystyle
K_{9} \leq\displaystyle
\varepsilon\int_0^t\left\|\nabla\cdot(\triangle{\bf
q}^\varepsilon)\right\|^2ds
+\varepsilon\int_0^t\left\|\nabla\cdot(\triangle\psi^\varepsilon)\right\|^2ds
\end{array}
\end{eqnarray}
and
\begin{eqnarray}\label{Le-conp-J17-c}
\arraycolsep=1.5pt
\begin{array}[b]{rl}\displaystyle
 K_{10}&\displaystyle\leq \int_0^t\left\|\nabla(\triangle\theta^\varepsilon)\right\|^2ds
 +C\int_0^t\left\|\triangle\psi^\varepsilon\widetilde{p}\right\|^2ds
 +C\int_0^t\left\|\nabla\psi^\varepsilon\cdot(\nabla\widetilde{p})\right\|^2ds
 \\[3mm]&\displaystyle+C\int_0^t\left\|\nabla^2\widetilde{p}\cdot\psi^\varepsilon\right\|^2ds
+C\int_0^t\left\|\triangle{\bf
q}\theta^\varepsilon\right\|^2ds \\[3mm]
&\displaystyle +C\int_0^t\left\|\nabla^2\theta^\varepsilon\cdot{\bf
q}\right\|^2ds +C\int_0^t\left\|\nabla{\bf
q}\cdot(\nabla\theta^\varepsilon)\right\|^2ds\\[3mm]
 &\displaystyle\leq \int_0^t\left\|\nabla(\triangle\theta^\varepsilon)\right\|^2ds
 +\sum\limits_{i=1}^{6}K_{10}^{i}.
\end{array}
\end{eqnarray}
Here $K_{10}^1-K_{10}^6$ are estimated as follows:
\begin{eqnarray}\label{Le-conp-2.83}
\arraycolsep=1.5pt
\begin{array}[b]{rl}\displaystyle
 K_{10}^1+K_{10}^2\leq &\displaystyle
 C\int_0^t\left(\left\|\widetilde{p}\right\|_{L^{\infty}}^2+\left\|\nabla\widetilde{p}\right\|_{L^{\infty}}^2\right)
 \left(\left\|\nabla\psi^\varepsilon\right\|^2+\left\|\triangle\psi^\varepsilon\right\|^2\right)ds\\[3mm]
\le& \displaystyle
C\int_0^t\left(\left\|\nabla\widetilde{p}\right\|_{H^1}^2+\left\|\nabla^2\widetilde{p}\right\|_{H^1}^2\right)
 \left(\left\|\nabla\psi^\varepsilon\right\|^2+\left\|\triangle\psi^\varepsilon\right\|^2\right)ds,
\end{array}
\end{eqnarray}

\begin{eqnarray}\label{Le-conp-2.84}
\arraycolsep=1.5pt
\begin{array}[b]{rl}\displaystyle
 K_{10}^3+K_{10}^4\leq &\displaystyle
 C\int_0^t\left(\left\|\psi^\varepsilon\right\|_{L^{\infty}}^2+\left\|\theta^\varepsilon\right\|_{L^{\infty}}^2\right)
 \left(\left\|\triangle\widetilde{p}\right\|^2+\left\|\triangle{\bf q}\right\|^2\right)ds\\[3mm]
 \le& \displaystyle C\int_0^t\left(\left\|\nabla\psi^\varepsilon\right\|_{H^{1}}^2+\left\|\nabla\theta^\varepsilon\right\|_{H^{1}}^2\right)
 \left(\left\|\triangle\widetilde{p}\right\|^2+\left\|\triangle{\bf q}\right\|^2\right)ds
\end{array}
\end{eqnarray}
and
\begin{eqnarray}\label{Le-conp-2.85}
\arraycolsep=1.5pt
\begin{array}[b]{rl}\displaystyle
 K_{10}^5+K_{10}^6\leq &\displaystyle
 C\int_0^t\left(\left\|{\bf q}\right\|_{L^{\infty}}^2+\left\|\nabla{\bf q}\right\|_{L^{\infty}}^2\right)
 \left(\left\|\nabla\theta^\varepsilon\right\|^2+\left\|\triangle\theta^\varepsilon\right\|^2\right)ds\\[3mm]
\le& \displaystyle
 C\int_0^t\left(\left\|\nabla{\bf q}\right\|_{H^{1}}^2+\left\|\nabla^2{\bf q}\right\|_{H^1}^2\right)
 \left(\left\|\nabla\theta^\varepsilon\right\|^2+\left\|\triangle\theta^\varepsilon\right\|^2\right)ds.
\end{array}
\end{eqnarray}
Plugging the estimates for $K_8\sim K_{10}$ into (\ref{Le-conp9-c}),
we get
\begin{eqnarray}\label{3Le-conpp-c}
\arraycolsep=1.5pt
\begin{array}[b]{rl}\displaystyle
 &\displaystyle\left\|\triangle\psi^\varepsilon\right\|^2+\left\|\triangle\theta^\varepsilon\right\|^2
+\int_0^t\left(\varepsilon\left\|\nabla\cdot(\triangle\psi^\varepsilon)\right\|^2
+\left\|\nabla^3\theta^\varepsilon\right\|^2\right)ds\\[3mm]
 \leq&\displaystyle C \varepsilon+C\int_0^t\left(\left\|\nabla\widetilde{p}\right\|^2_{H^2}+\left\|\nabla{\bf q}\right\|^2_{H^2}\right)
\left(\left\|\psi^\varepsilon\right\|_{H^2}^2+\left\|\theta^\varepsilon\right\|_{H^2}^2\right)ds.
\end{array}
\end{eqnarray}
Combination of (\ref{Le-conp8-4.10}), (\ref{3Le-conpr-cc}) and
(\ref{3Le-conpp-c}) yields
\begin{eqnarray}\label{3Le-conpf-c}
\arraycolsep=1.5pt
\begin{array}[b]{rl}\displaystyle
 &\displaystyle\left\|\psi^\varepsilon\right\|_{H^2}^2+\left\|\theta^\varepsilon\right\|_{H^2}^2
+\int_0^t\left(\varepsilon\left\|\nabla\psi^\varepsilon\right\|_{H^2}^2
+\left\|\nabla\theta^\varepsilon\right\|_{H^2}^2\right)ds\\[3mm]
 \leq&\displaystyle C \varepsilon+C\int_0^t\left(\left\|\nabla\widetilde{p}\right\|^2_{H^2}+\left\|\nabla{\bf q}\right\|^2_{H^2}\right)
\left(\left\|\psi^\varepsilon\right\|_{H^2}^2+\left\|\theta^\varepsilon\right\|_{H^2}^2\right)ds.
\end{array}
\end{eqnarray}
which, together with Theorem \ref{theorem 1 } and Gronwall¡¯s
inequality gives (\ref{Le-con1lemma-c}). This completes the proof of
Lemma \ref{Le-3.6}. \endpf

\bigskip

From Lemma \ref{Le-3.6},  we get (\ref{Th-con-cb}). Using Sobolev
inequality and (\ref{Th-con-cb}), we get (\ref{theorem 2.2 last}).
This completes the proof of Theorem \ref{theorem 2 }.

\section*{Acknowledgments}\ \ The research was
supported by the National Natural Science Foundation of China
$\#$10625105, $\#$11071093, the PhD specialized grant of the
Ministry of Education of China $\#$20100144110001, and the Special
Fund for Basic Scientific Research  of Central Colleges
$\#$CCNU10C01001, CCNU12C01001.

\end{document}